\documentclass[final,hidelinks,onefignum,onetabnum]{siamart250211}

\usepackage{amsfonts,amssymb}
\usepackage{graphicx}
\ifpdf
  \DeclareGraphicsExtensions{.eps,.pdf,.png,.jpg}
\else
  \DeclareGraphicsExtensions{.eps}
\fi
\usepackage{bm}
\usepackage{tikz}
\usetikzlibrary{math,positioning}

\newsiamremark{remark}{Remark}
\newsiamthm{stdass}{Standard Assumptions}
\newsiamthm{conjecture}{Conjecture}

\newcommand{\eps}{\epsilon}
\newcommand{\RR}{\mathbb{R}}

\newcommand{\grad}{\nabla}
\newcommand{\Div}{\nabla\cdot}

\newcommand{\bg}{\mathbf{g}}
\newcommand{\bn}{\mathbf{n}}
\newcommand{\bq}{\mathbf{q}}
\newcommand{\bu}{\mathbf{u}}
\newcommand{\bv}{\mathbf{v}}

\newcommand{\bU}{\mathbf{U}}

\newcommand{\bzero}{\bm{0}}

\newcommand{\cK}{\mathcal{K}}
\newcommand{\cQ}{\mathcal{Q}}
\newcommand{\cT}{\mathcal{T}}
\newcommand{\cV}{\mathcal{V}}
\newcommand{\cX}{\mathcal{X}}

\newcommand{\nn}{{\text{\textnormal{n}}}}
\newcommand{\pp}{{\text{\textnormal{p}}}}
\newcommand{\qq}{{\text{\textnormal{q}}}}
\newcommand{\rr}{{\text{\textnormal{r}}}}

\newcommand{\rhoi}{\rho_{\text{i}}}

\newcommand{\nsubset}{\not\subset}

\ifpdf
\hypersetup{
  pdftitle={Surface elevation errors in finite element Stokes models for glacier evolution},
  pdfauthor={E. Bueler}
}
\fi

\headers{Surface elevation errors in glacier models}{E. Bueler}

\title{Surface elevation errors in finite element Stokes models for glacier evolution}

\author{Ed Bueler\thanks{Dept.~Mathematics \& Statistics, University of Alaska Fairbanks, USA (\email{elbueler@alaska.edu}).}}

\begin{document}

\maketitle

\begin{abstract}
The primary data which determine the evolution of glaciation are the bedrock elevation and the surface mass balance.  From this data, which we assume is defined over a fixed land region, the glacier's geometry solves a free boundary problem which balances the time derivative of the surface elevation, the surface velocity from the Stokes flow of the ice, and the surface mass balance rate.  This problem can be posed in weak form as a variational inequality over a cone of admissible surface elevation functions, those which are above the bedrock topography.  After some preparatory theory for the Stokes problem, we conjecture that the corresponding continuous-space, implicit time-step variational inequality problem is well-posed if the surface kinematical equation is appropriately regularized.  This conjecture is supported by physical arguments and numerical evidence.  We then prove a general theorem which bounds the numerical error made by finite element approximations of nonlinear-operator variational inequalities in Banach spaces.  This bound is a sum of error terms of different types, special to variational inequalities.  When it is applied to the implicit time-step glacier problem there are three terms in the bound: an error from discretizing the bed elevation primarily near the glacier margin, an error from numerically solving for the Stokes velocity, and finally a C\'ea-type best-approximation error for the finite element space representation of the surface elevation.  The design of glacier models is then reconsidered based on this mostly \emph{a priori} error analysis.
\end{abstract}



\section{Introduction} \label{sec:intro}

Glacier and ice sheet simulations model the flowing ice as a free-surface layer of very-viscous, incompressible, and non-Newtonian fluid \cite{GreveBlatter2009,SchoofHewitt2013}.  The two essential input data into such simulations are the bedrock elevation, which is assumed for simplicity to be independent of time, and the time- and space-dependent surface mass balance rate (SMB), the climate.  By definition, the SMB is the balance between accumulating snow and the loss of melt water, through runoff, at the upper surface of the glacier \cite{Cogleyetal2011}.  Note that elevations are measured in meters, and SMB is in (ice-equivalent) meters per second.

Time-dependent glacier simulations \cite{SchoofHewitt2013} produce the glacier's evolving geometry and flow velocity, though comprehensive models \cite[for example]{Winkelmannetal2011} have additional physics and outputs.  For example, they might track the enthalpy/temperature \cite{Aschwandenetal2012} of the ice, or liquid water along glacier surfaces.  However, the current work only addresses conservation of mass and momentum, but not energy.  Liquid water plays no role.  We only consider simulations of land-based glaciers and ice sheets (continent-scale glaciers), without floating portions. Furthermore we will assume zero velocity at the base of the ice, a non-sliding and non-penetrating condition.  On the other hand, we will not make any of the shallowness assumptions which are common in present-day comprehensive models, and furthermore we will address the mathematical model rigorously.

One may parameterize the glacier's geometry using either the (upper) surface elevation or the ice thickness.  At any time and map-plane location where a glacier exists the surface elevation exceeds the elevation of the bedrock, equivalently the ice thickness is positive.  The computed flow velocity is only defined at those locations and times where ice is present, on an evolving 3D domain between the bedrock and surface elevations.  In other words, admissibility of the surface elevation is required to make the velocity and pressure Stokes problem meaningful.

\begin{figure}[ht]
\centering
\includegraphics[width=0.6\textwidth]{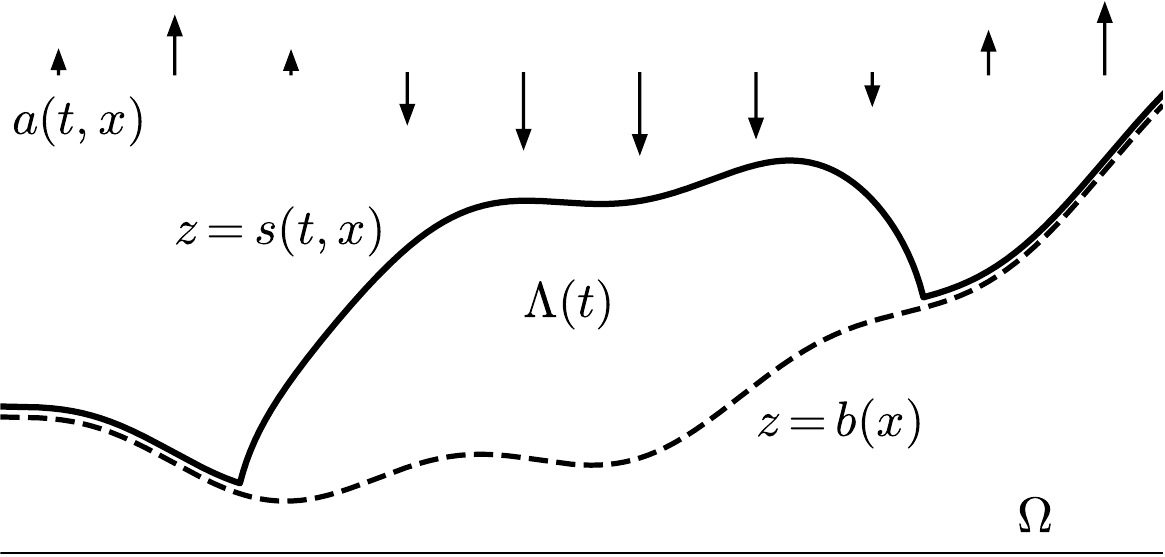}
\caption{Glacier notation: $x\in\Omega\subset\RR^2$ and $(x,z) \in \Lambda(t)\subset\RR^3$.}
\label{fig:stokesdomain}
\end{figure}

Our notation is sketched in Figure \ref{fig:stokesdomain}.  Let $\Omega \subset \RR^2$ be a fixed domain of land, with piecewise-$C^1$ boundary $\partial\Omega$, and with map-plane coordinates $x=(x_1,x_2)\in\Omega$.  We are given, as data, a continuous bedrock elevation function $b(x)$ on $\Omega$, and a continuous SMB function $a(t,x)$ on $[0,T]\times \Omega$, for some $T>0$.  Where $a>0$ (accumulation; downward arrows in Figure \ref{fig:stokesdomain}), a glacier will exist.  However, where $a<0$ (ablation; upward arrows) then either ice exists with an ablating surface, supported by flow from accumulation areas, or no glacier exists.  Determining which situation applies at a given location and time, and the surface elevation there, requires solving a free-boundary problem.

Let $s(t,x)$ be the (solution) ice surface elevation.  This is defined for all $x\in\Omega$, and it is subject to the constraint that the surface $z=s$ must be at or above the bedrock ($s \ge b$).  In regions with no ice, $s=b$ holds.

The solution ice velocity $\bu(t,x,z)$ and pressure $p(t,x,z)$ are defined only on the open 3D domain (Figure \ref{fig:stokesdomain})
\begin{equation}
\Lambda(t) = \left\{(x,z)\,:\,b(x) < z < s(t,x)\right\} \subset \Omega \times \RR. \label{eq:icydomain}
\end{equation}
This aspect of glacier modeling deserves emphasis:  The time-dependent 3D domain $\Lambda(t)$, on which the velocity and pressure are meaningful, is determined by the evolving surface elevation $s$, which is itself part of the coupled model solution.  The surface trace of the ice velocity---a Sobolev space for this field appears in Section \ref{sec:stokes}---is extended by zero, so that it is defined for all $(t,x) \in [0,T]\times\Omega$:
\begin{equation}
\bu|_s(t,x) = \begin{cases} \bu(t,x,s(t,x)), & s(t,x)>b(x) \\
                            \bzero, & \text{otherwise} .\end{cases} \label{eq:defineus}
\end{equation}
Compare flux extension by zero in \cite{SchoofHewitt2013}.

Now the glacier geometry evolution model says that an infinite-dimensional nonlinear complementarity problem (NCP) \cite{FacchineiPang2003}, an obstacle problem, applies in $[0,T]\times \Omega$.  Defining $\bn_s = \left<-\partial s/\partial x_1,-\partial s/\partial x_2,1\right>$, an upward surface normal vector, the NCP is
\begin{subequations}
\label{eq:ncp}
\begin{align}
s - b &\ge 0 \label{eq:ncp:constraint} \\
\frac{\partial s}{\partial t} - \bu|_s \cdot \bn_s - a &\ge 0 \label{eq:ncp:residualpos} \\
(s - b) \left(\frac{\partial s}{\partial t} - \bu|_s \cdot \bn_s - a\right) &= 0.
\end{align}
\end{subequations}
Using \eqref{eq:defineus}, this system implies that either a location is ice free ($s=b$), where the climate is locally ablating ($a\le 0$), or the surface kinematical equation holds:
\begin{equation}
\frac{\partial s}{\partial t} - \bu|_s \cdot \bn_s - a = 0.  \label{eq:ske}
\end{equation}
This equation says that the (non-material) ice surface moves vertically according to the sum of an ice velocity component at the surface and the SMB.  This statement of mass conservation \cite{Aschwandenetal2012}, also called the free-surface equation \cite{LofgrenAhlkronaHelanow2022}, is a standard description of glacier geometry evolution in numerical models \cite{GreveBlatter2009,SchoofHewitt2013}.  Though it is not yet common to state system \eqref{eq:ncp} completely as an NCP, as here, glaciologists would agree with its conditions.  While condition \eqref{eq:ncp:constraint} is sometimes stated in glacier literature \cite{Durandetal2009,Halfar1981,JouvetBueler2012,PiersantiTemam2023,WirbelJarosch2020}, the complementary fact \eqref{eq:ncp:residualpos}, that the residual is everywhere nonnegative, though observed already in \cite{Calvoetal2003}, is rarely written \cite{SchoofHewitt2013}.

Because a simulated glacier needs to be able to advance into unglaciated locations, we assume that the SMB $a$ is defined everywhere in $\Omega$, regardless of whether a glacier is present or not.  In ice-free areas $a$ should have the value which a glacier surface at elevation $b$ would experience.  For example, it can be computed, from precipitation and an energy balance model \cite{GreveBlatter2009}, by balancing snow accumulation minus the ablation caused by using the available energy for melt if ice were present.

The non-shallow ice dynamics model used in this paper is the non-sliding (e.g.~frozen base), isothermal, non-Newtonian, and incompressible Stokes system \cite{GreveBlatter2009,JouvetRappaz2011,SchoofHewitt2013}.  It is applied over the domain $\Lambda(t)$ defined in \eqref{eq:icydomain}.  Let $\Gamma_s(t)$, $\Gamma_b(t)$ denote the upper and lower surfaces at $z=s$ and $z=b$, respectively.  As $s$ evolves in time the upper surface moves vertically, but also both boundary \emph{sets} $\Gamma_s(t)$, $\Gamma_b(t)$ evolve in time as the lateral extent of the ice changes.  The possibility of cliffs at the ice margin will be neglected (Section \ref{sec:model}), so $\partial \Lambda(t) = \overline{\Gamma_s(t)} \cup \overline{\Gamma_b(t)}$.

Now let $D\bu=(\grad \bu + \grad \bu^{\top})/2$ denote the strain rate tensor, with Frobenius norm $|D\bu| = \left((D\bu)_{ij} (D\bu)_{ij}\right)^{1/2}$.  A shear-thinning constitutive relation, a regularized Glen's flow law \cite{GreveBlatter2009}, computes the effective dynamic ice viscosity:
\begin{equation}
\nu(D\bu) = \nu_\pp \left(|D\bu|^2 + \mu_0\right)^{(\pp-2)/2}. \label{eq:glen}
\end{equation}
The exponent $1 < \pp \le 2$, written $\pp=(1/\nn)+1$ in terms of Glen's exponent $\nn\ge 1$, is commonly taken to be $\pp=4/3$ \cite{GreveBlatter2009}.  The coefficient $\nu_\pp>0$ necessarily has $\pp$-dependent units, but $\nu(D\bu)$ has SI units $\text{kg}\,\text{m}^{-1}\,\text{s}^{-1}$.  The values of $\pp$ and $\nu_\pp$, assumed to be constant here, are constrained by laboratory experiments \cite{GoldsbyKohlstedt2001,GreveBlatter2009}.  While $\pp=2$ yields a Newtonian fluid with constant viscosity, when $\pp < 2$ the $\mu_0>0$ regularization implies that $\nu(D\bu)$ is bounded.

In the Stokes model the velocity and pressure solve the following equations:%
\begin{subequations}
\label{eq:stokes}
\begin{align}
- \nabla \cdot \left(2 \nu(D\bu)\, D\bu\right) + \nabla p &= \rhoi \bg && \text{stress balance in $\Lambda(t)$} \\
\nabla \cdot \bu &= 0 && \text{incompressible in $\Lambda(t)$} \label{eq:stokes:incomp} \\
\left(2 \nu(D\bu) D\bu - pI\right) \bn_s &= \bzero && \text{stress-free on $\Gamma_s(t)$}\label{eq:stokes:stressfreesurface} \\
\bu  &= \bzero && \text{non-sliding on $\Gamma_b(t)$} \label{eq:stokes:noslide}
\end{align}
\end{subequations}
The density of ice $\rhoi$ and the acceleration of gravity $\bg$ are here assumed constant.  This Stokes problem \eqref{eq:glen}, \eqref{eq:stokes} is commonly discretized using finite elements (FEs) \cite{IsaacStadlerGhattas2015,Jouvetetal2008,Pattynetal2008}, and we will return to its FE approximation in Section \ref{sec:application}.

In summary at this point, this work considers an evolving free-surface flow for a glacier, subject to a signed climate which can add or remove ice, simultaneously with a non-Newtonian Stokes problem which must be solved within the evolving 3D domain of ice.  This coupled initial-boundary value problem, in strong form as \eqref{eq:icydomain}--\eqref{eq:stokes}, requires data $b(x)$ and $a(t,x)$, plus an initial surface elevation $s(0,x)$.  The solution variables are $s(t,x)$, which is defined everywhere over $[0,T]\times \Omega$ and subject to $s \ge b$, and $\bu(t,x,z)$ and $p(t,x,z)$, which are defined on $\Lambda(t)$ given by \eqref{eq:icydomain}.

\begin{figure}
\begin{minipage}[t]{0.8\textwidth}
\vspace{0pt}
\includegraphics[width=\textwidth]{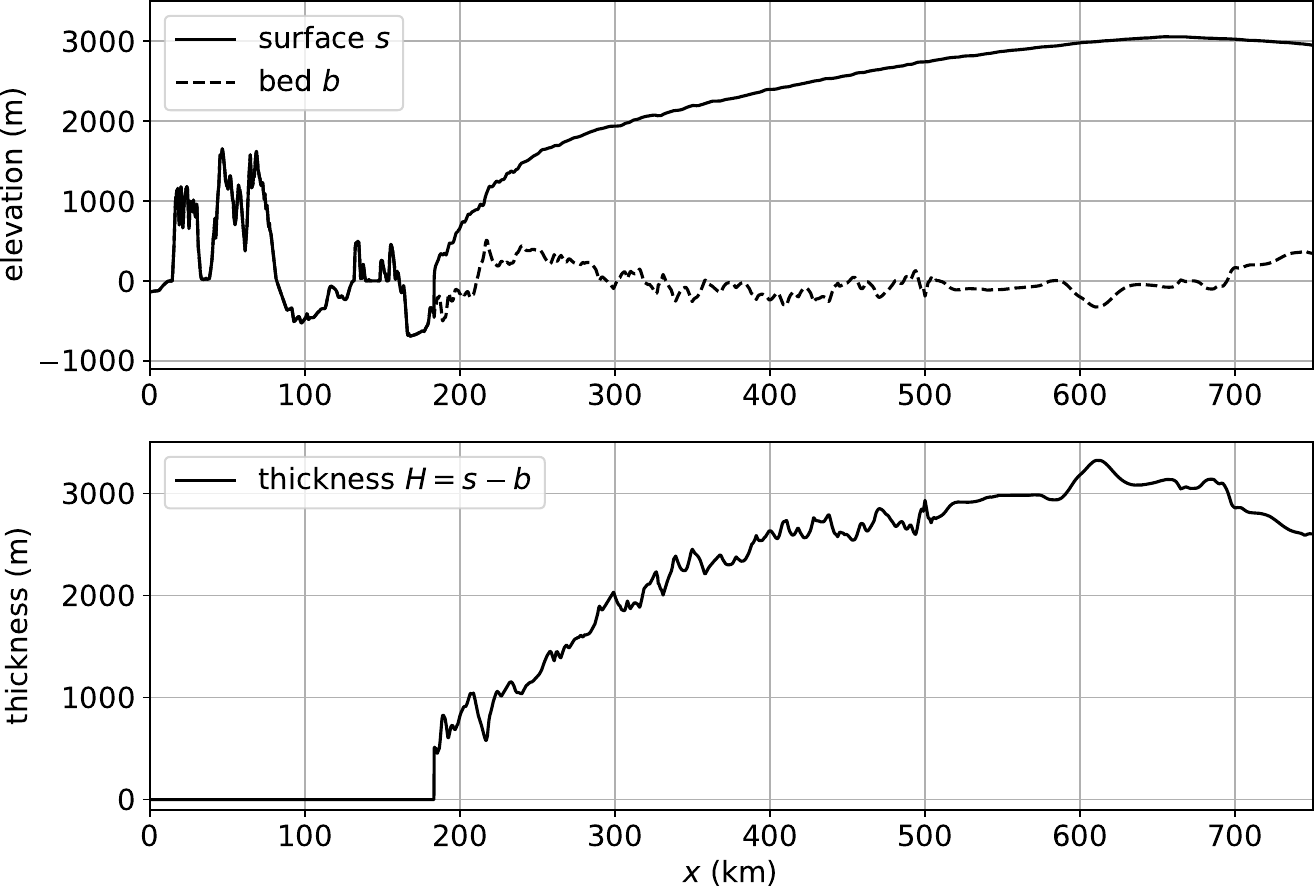}
\end{minipage}
\,
\begin{minipage}[t]{0.15\textwidth}
\vspace{10pt}
\includegraphics[width=\textwidth]{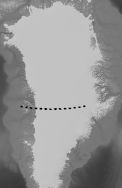}
\end{minipage}
\caption{A cross-section of the Greenland ice sheet at $70^\circ$N latitude \cite{Morlighemetal2017}; see inset.  Top: While the ice surface $s$ is smoothed because of ice flow, the bedrock elevation $b$ is much rougher.  Bottom: The corresponding ice thickness $H = s-b$ ``inherits'' the low regularity of $b$.}
\label{fig:giscross}
\end{figure}

A formulation using ice thickness to parameterize glacier geometry would have a different character from ours.  Surface elevation is subject to a flow-caused smoothing effect, illustrated for a real ice sheet in Figure \ref{fig:giscross}.  For land-based glaciers $s(t,x)$ is (spatially) smoother than the thickness $H(t,x) = s(t,x)-b(x)$; the latter ``inherits'' the lower regularity of the steep and eroded bedrock topography.

Observe that NCP \eqref{eq:ncp} is the only part of the model where a time derivative appears.  Indeed, because the flow is very-viscous \cite{Acheson1990}, the Stokes sub-model acts as an instantaneous ``algebraic'' constraint on the evolution.  Therefore the coupled and infinite-dimensional glacier evolution problem is simultaneously a differential algebraic equation (DAE) system \cite{AscherPetzold1998,LofgrenAhlkronaHelanow2022} and a free-boundary, obstacle-type NCP.  While this is a now-standard physical model, for an evolving glacier, as a mathematical model it is poorly understood.  No well-posedness theory is known to this author, although existence is known for a shallow limit \cite{JouvetBueler2012,PiersantiTemam2023}, with non-uniqueness in general (Appendix \ref{app:noncoercive}).  The current paper seems to be the first where the non-shallow model is stated with (proposed) function spaces for all variables, including the surface elevation.

Furthermore, to the author's knowledge all existing evolution models using Stokes dynamics, or with any other form of a longitudinal stress balance, use an explicit or semi-implicit time-stepping scheme for the geometry \cite[for examples]{Brinkerhoff2023, Bueler2023, Durandetal2009, Jouvetetal2008, LofgrenAhlkronaHelanow2022, WirbelJarosch2020}, with a recent exception which is fully-implicit on a fixed map-plane region \cite{AhlkronaLofgrenHenry2025}.  Specifically, when implicitness is applied in these implementations the portion of $\RR^2$ covered by (solution) ice is not general.  Many schemes limit margin advance or retreat to one cell from the previous time step, though there is no physical justification for such a restriction.  The current paper considers a single time step of glacier evolution, without restrictions on the new (solution) ice-covered area.

For finite-dimensional DAE systems, implicit schemes are the standard choice because they handle unbounded stiffness \cite{AscherPetzold1998}.  In our infinite-dimensional case, where NCP \eqref{eq:ncp} is constrained by the ``algebraic'' Stokes problem \eqref{eq:glen}--\eqref{eq:stokes}, implicit schemes are also the natural choice.  We consider a backward Euler semi-discretization for simplicity.  For a time step $\Delta t > 0$, and times $\{t_n\}$ with $t_n-t_{n-1}=\Delta t$, this scheme applied to \eqref{eq:ncp} yields a new NCP for the updated surface elevation $s^n \approx s(t_n,x)$:%
\begin{subequations}
\label{eq:be:ncp}
\begin{align}
s^n - b &\ge 0 \label{eq:be:ncp:constraint} \\
s^n - \Delta t\,\bu|_{s^n} \cdot \bn_{s^n} - \ell^n &\ge 0 \label{eq:be:ncp:residualpos} \\
(s^n - b) \left(s^n - \Delta t\,\bu|_{s^n} \cdot \bn_{s^n} - \ell^n\right) &= 0. \label{eq:be:ncp:complementarity}
\end{align}
\end{subequations}
For clarity we have collected a source term\footnote{Quadrature could be used on the integral of SMB, e.g.~$\int_{t_{n-1}}^{t_n} a(t,x)\,dt \approx \Delta t\,a(t_n,x)$, but the integral form is preferred because of fast temporal variation in $a$; compare \cite{Cogleyetal2011}.} $\ell^n(x) = s^{n-1}(x) + \int_{t_{n-1}}^{t_n} a(t,x)\,dt$.  In Section \ref{sec:model} we will re-write the coupled NCP problem \eqref{eq:glen}--\eqref{eq:be:ncp} as a weak-form variational inequality (VI) for $s^n$ in an admissible subset of a Banach space.

Backward Euler scheme \eqref{eq:be:ncp} is merely the simplest implicit choice.  Extensions to higher-order A-stable and stiff decay methods \cite{AscherPetzold1998} are natural and straightforward.  In the other direction, a (fully) explicit version of \eqref{eq:be:ncp} \cite{Lengetal2012} replaces $s^n$ in the surface motion term by the previous surface elevation: $\bu|_{s^n} \cdot \bn_{s^n} \to \bu|_{s^{n-1}} \cdot \bn_{s^{n-1}}$.  One possible semi-implicit scheme \cite{Durandetal2009} uses the surface velocity from the old time, but with the updated value for the surface slope part: $\bu|_{s^n} \cdot \bn_{s^n} \to \bu|_{s^{n-1}} \cdot \bn_{s^n}$.  A more stable form of semi-implicitness comes from modifying the body force in the coupled Stokes problem using an explicit surface estimate \cite{LofgrenAhlkronaHelanow2022}.

Regarding simulation performance, one can make a case for implicit time-stepping based on a computational complexity analysis \cite{Bueler2023}.  Whether it leads to actual performance improvements will depend on positive answers to key questions:  Is the implicit-step VI problem well-posed?  Can it be solved accurately and efficiently?  Does the (coupled) implicit scheme have unconditional stability?  In the current paper we conjecture well-posedness and address FE accuracy.  Time-stepping stability, time-dependent convergence, and solver efficiency are all topics for further research.

This paper is organized as follows.  Section \ref{sec:stokes} recalls the theory of the power-law Stokes problem, on a fixed 3D domain, for which we give a new bound on the surface trace of the velocity solution (Lemma \ref{lem:surfacetracebound}).  In Section \ref{sec:model} we re-formulate the coupled, implicit-step NCP problem \eqref{eq:glen}--\eqref{eq:be:ncp} as a VI weak form.  The key coupling is the surface motion term in \eqref{eq:be:ncp}, so in Section \ref{sec:conjectural} we hypothesize the coercivity of a regularized version (Conjecture \ref{conj:regcoercive}).  Physical and modeling experience, along with new numerical results (Appendix \ref{app:numerical}), support this conjecture.  Section \ref{sec:abstractestimate} considers the FE approximation of abstract VI problems in Banach spaces, a class which includes the temporally-discrete, but spatially-continuous, regularized VI problem \eqref{eq:regularizedvi}.  An error estimate for FE solutions of VI problems is proven (Theorem \ref{thm:abstractestimate}).  Though it makes coercivity and Lipschitz assumptions, this apparently-new bound permits nonlinear operators on Banach spaces, so it significantly extends a classical bilinear result by Falk \cite{Falk1974}.  In Section \ref{sec:application} we apply Theorem \ref{thm:abstractestimate} to the glacier problem, yielding our main result, Theorem \ref{thm:glacierapp}, which is a bound on FE error when solving \eqref{eq:regularizedvi}.  Meanings for each term in the bound, and consequences for glacier model and FE method choices, are addressed in Sections \ref{sec:application} and \ref{sec:conclusion}.  Note that we will use only seven abbreviations which are standard in their fields: DAE (differential-algebraic equations), FE (finite element), NCP (nonlinear complementarity problem), PDE (partial differential equation), SIA (shallow ice approximation), SMB (surface mass balance), and VI (variational inequality).

\section{Surface velocity from the Stokes sub-model} \label{sec:stokes}

In this Section we address only the Stokes sub-model \eqref{eq:glen}--\eqref{eq:stokes}, applied on a fixed domain $\Lambda \subset \RR^3$, e.g.~defined by \eqref{eq:icydomain} at a fixed time $t$.  This sub-model computes the surface velocity $\bu|_s$ in NCP \eqref{eq:be:ncp}, so we return to the time-dependent coupling in the next section.

Suitable function spaces for well-posedness of this sub-model are known.  Recall that $1 < \pp \le 2$.  Denote the Sobolev space \cite{Evans2010} of real-valued functions on $\Lambda$ with $\pp$th-power integrable first derivatives by $W^{1,\pp}(\Lambda)$.  Let $\cV = W_0^{1,\pp}(\Lambda; \RR^3)$ be the corresponding vector-valued functions with trace zero along the ice base $\Gamma_b\subset\partial \Lambda$, which is assumed to have positive measure.  For $[H]\ge 1$ a representative vertical glacier scale in meters, we define the norm
\begin{equation}
\|\bv\|_{\cV} = \left(\int_\Lambda |\bv|^\pp\,dx\,dz + [H]^\pp \int_\Lambda |\grad\bv|^\pp\,dx\,dz\right)^{1/\pp}. \label{eq:vnorm}
\end{equation}
Here $|\grad\bv|=\left(\grad\bv : \grad\bv\right)^{1/2}$ is the Frobenius norm on $\RR^{3\times 3}$, for $A:B=a_{ij}b_{ij}$.  The scale $[H]$ gives $\|\bv\|_{\cV}$ consistent units; see Remark 1.2.1 in \cite{BoffiBrezziFortin2013}.  The volume element $dx\,dz = dx_1\,dx_2\,dz$ will be suppressed from now on.

Let $\cQ=L^{\pp'}(\Lambda)$ where $\pp'=\pp/(\pp-1)$.  For $(\bu,p), (\bv,q) \in \mathcal{M} = \cV \times \cQ$, the mixed velocity-pressure space, define
\begin{equation}
F_\Lambda(\bu,p)[\bv,q] = \int_\Lambda 2 \nu(D\bu) D\bu : D\bv - p \Div\bv - (\Div\bu) q - \rhoi \bg \cdot \bv. \label{eq:glenstokes:fcnl}
\end{equation}
The weak form of the Stokes sub-model \eqref{eq:glen}--\eqref{eq:stokes} now seeks $(\bu,p) \in \mathcal{M}$ satisfying
\begin{equation}
F_\Lambda(\bu,p)[\bv,q] = 0 \qquad \text{for all } (\bv,q) \in \mathcal{M}. \label{eq:glenstokes:weak}
\end{equation}

Problem \eqref{eq:glenstokes:weak} is well-posed if $\partial\Lambda$ is appropriately regular, e.g.~piecewise $C^1$ \cite{JouvetRappaz2011} or polygonal \cite{Belenkietal2012}.  We assume throughout that $\Lambda$ is a bounded Lipschitz domain.  A well-posedness proof uses a specific constitutive relation; our regularization \eqref{eq:glen} of the Glen law  matches that in \cite{Belenkietal2012,IsaacStadlerGhattas2015}.  There exists a unique pair $(\bu,p) \in \mathcal{M}$ solving \eqref{eq:glenstokes:weak}, with $\bu\in \cV_0 = \{\bv\in\cV\,:\, \Div\bv=0\}$.  The following \emph{a priori} bound is proven, for completeness, in Appendix \ref{app:provestokesapriori}.

\begin{lemma} \label{lem:stokesapriori}
There is $C>0$ depending on $\pp$, $\rhoi |\bg|$, $\nu_\pp$, $\mu_0$, $[H]$, and $\Lambda$ so that if $\bu\in\cV_0$ is the solution from \eqref{eq:glenstokes:weak} then
\begin{equation}
\|\bu\|_{\cV} \le C. \label{eq:stokesapriori}
\end{equation}
\end{lemma}

Our primary purpose is to study NCP \eqref{eq:be:ncp} for glacier geometry, specifically its weak form (next section), and for that we need control on the surface trace $\bu|_s$.  Because $\Lambda$ is a bounded Lipschitz domain, the trace theorem \cite{Grisvard1985} gives the following Lemma; compare \cite[Theorem 5.5.1]{Evans2010} for $C^1$ boundaries.

\begin{lemma}[Trace inequality] \label{lem:trace}
There exists a constant $C$, depending on $\pp$, $[H]$, and $\Lambda$, so that for all $\bv \in \cV$,
\begin{equation}
\int_{\Gamma_s} \big|\bv|_s\big|^\pp \,dS \le C \|\bv\|_{\cV}^\pp. \label{eq:trace}
\end{equation}
Here $dS$ denotes the area element on $\partial\Lambda$.
\end{lemma}

Combining Lemmas \ref{lem:stokesapriori} and \ref{lem:trace} yields the following bound on the surface trace of the velocity.  When applying this result in the next section, recall that $\Lambda=\Lambda(t)$ and $\Gamma_s=\Gamma_s(t)$ are defined via \eqref{eq:icydomain}, in terms of $s(t,x)$ and $b(x)$.

\begin{lemma}[Surface velocity bound] \label{lem:surfacetracebound}
There is a constant $C>0$, computable from physical constants and the geometry of $\Lambda$, so that for $\bu\in\cV_0$ solving \eqref{eq:glenstokes:weak},
\begin{equation}
\int_{\Gamma_s} \big|\bu|_s\big|^\pp \,dS \le C. \label{eq:surfacetracebound}
\end{equation}
\end{lemma}

\section{The weak form of an implicit time step} \label{sec:model}

Suppose $(t_n)$ are increasing times in $[0,T]$, with $t_0=0$, and write $\Delta t = t_n-t_{n-1}$ for a generic step.  Let $a^n(x)$ be the temporal average of the climatic data $a(t,x)$ over $[t_{n-1},t_n]$.  Suppose that $s^n(x)\approx s(t_n,x)$ approximates the surface elevation at time $t_n$.  The backward Euler scheme \cite{AscherPetzold1998} for the surface kinematical equation \eqref{eq:ske} is
\begin{equation}
\frac{s^n - s^{n-1}}{\Delta t} - \bu|_{s^n} \cdot \bn_{s^n} - a^n = 0. \label{eq:be:ske}
\end{equation}
Note how $s^n$ appears in both the surface velocity $\bu|_{s^n}$ and the slope $\bn_{s^n}$.

For cleaner appearance in what follows we will write $s=s^n$ for the unknown surface elevation.  As noted in the Introduction, $s$ solves a problem of free-boundary type, namely NCP \eqref{eq:be:ncp} in strong form, which says that either there is bare ground ($s=b$) or equation \eqref{eq:be:ske} holds.  Leaving the precise Banach space $\cX$ to be determined, admissible surface elevations for this problem form a convex and closed cone:
\begin{equation}
\boxed{\cK = \left\{\sigma \in\cX\,:\,\sigma|_{\partial\Omega}=b|_{\partial\Omega} \text{ and } \sigma \ge b\right\}.}  \label{eq:be:admissible}
\end{equation}
(The Dirichlet boundary condition is included in the definition of $\cK$.)  We assume throughout that $b\in C^1(\bar\Omega)$.

The weak form of the free-boundary problem for \eqref{eq:be:ske} is derived by assuming that $s \in \cK$ is a sufficiently-regular solution of \eqref{eq:be:ncp}.  For tidiness collect a source term defined over all of $\Omega$:
\begin{equation}
\boxed{\ell^n(x) = s^{n-1}(x)+\Delta t\,a^n(x).} \label{eq:be:source}
\end{equation}
Let $\Omega_I = \{x\in\Omega\,:\,s(x)>b(x)\}$ be the (measurable) subset on which the constraint \eqref{eq:be:ncp:constraint} is inactive, i.e.~where glacier ice is present.  By \eqref{eq:be:ncp:complementarity}, integration over $\Omega_I$ gives%
\begin{equation}
\int_{\Omega_I} \left(s - \Delta t\,\bu|_s \cdot \bn_s - \ell^n\right)\,(\sigma-s) = 0  \label{eq:inactivetruth}
\end{equation}
for any $\sigma\in\cK$.  On the other hand, let $\Omega_A = \Omega \setminus \Omega_I$ be the active (ice-free) region.  Using extension by zero \eqref{eq:defineus}, inequality \eqref{eq:be:ncp:residualpos} over $\Omega_A$ implies $b-\ell^n = s - \Delta t\,\bu|_s \cdot \bn_s - \ell^n \ge 0$.  Since also $\sigma-s=\sigma-b\ge 0$ on $\Omega_A$, integration yields an inequality:
\begin{equation}
\int_{\Omega_A} \left(s - \Delta t\,\bu|_s \cdot \bn_s - \ell^n\right)\,(\sigma-s) \ge 0.  \label{eq:activetruth}
\end{equation}
Adding \eqref{eq:inactivetruth} and \eqref{eq:activetruth} gives a variational inequality (VI; \cite{KinderlehrerStampacchia1980})
\begin{equation}
\int_\Omega \left(s - \Delta t\,\bu|_s \cdot \bn_s - \ell^n\right)\,(\sigma-s) \ge 0 \quad \text{for all } \sigma \in \cK. \label{eq:be:viearly}
\end{equation}
This VI holds true for $s\in\cK$, in advance of knowing which part of $\Omega$ is ice-covered.

\begin{figure}[ht]
\begin{center}
\begin{tikzpicture}[scale=1.0]
  \node (wedge) {\includegraphics[width=25mm]{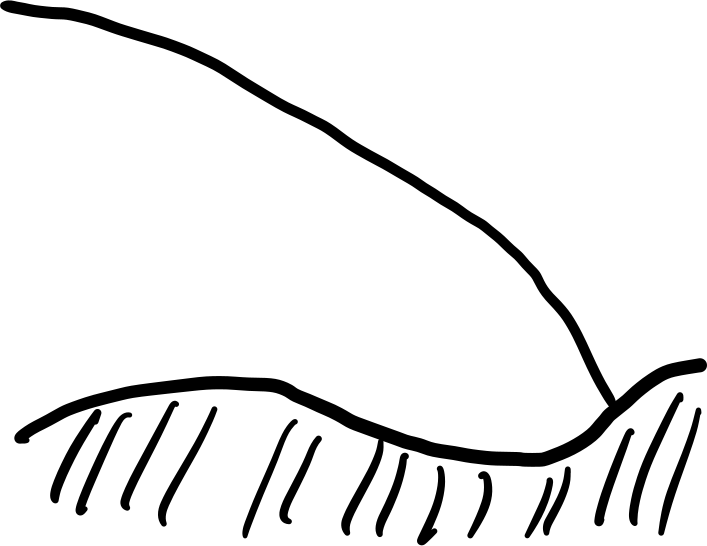}} node[xshift=-4mm, yshift=1mm] at (wedge.center) {{\small \emph{ice}}};
  \node[right=of wedge] (unbounded) {\includegraphics[width=26mm]{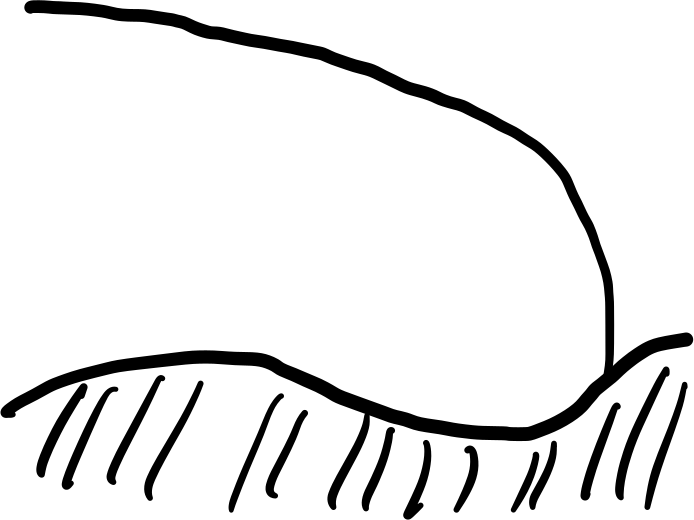}} node[xshift=-2mm, yshift=1mm] at (unbounded.center) {{\small \emph{ice}}};
  \node[right=of unbounded] (realistic) {\includegraphics[width=26mm]{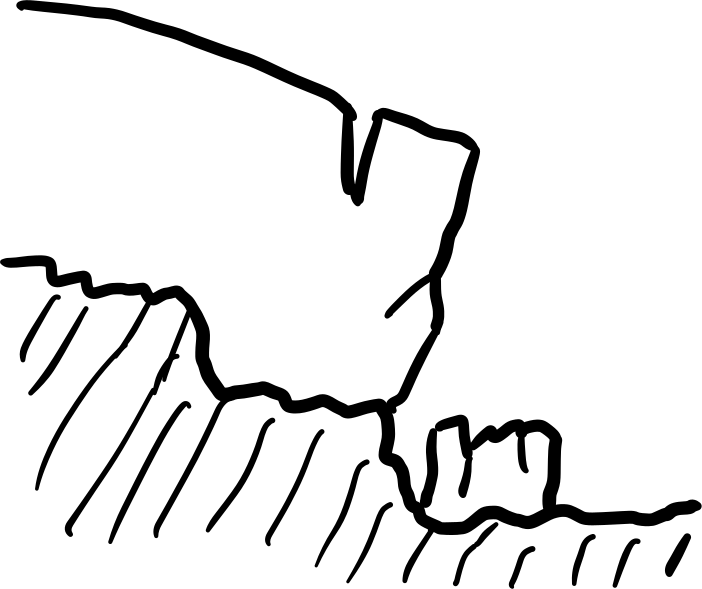}} node[xshift=-6mm, yshift=4mm] at (realistic.center) {{\small \emph{ice}}};
\end{tikzpicture}
\end{center}

\vspace{-2mm}

\caption{Glacier margins with finite-slope ``wedge'' shapes (left) are possible, while shallow theory has fractional powers (center).  Actual glacier margins might have overhangs and fractures (right).}
\label{fig:margins}
\end{figure}

It remains to identify a Sobolev space $\cX$ suitable for surface elevations.  However, within a physical Stokes-based theory, the shape that should be predicted for a glacier's grounded margin is not clear (Figure \ref{fig:margins}).  ``Wedge'' shapes with bounded gradients have been hypothesized from observations \cite{EchelmeyerKamb1986}, while shallow-ice theory suggests root-type (fractional-power) shapes, with different powers for advance and retreat \cite{Bueleretal2005,JouvetBueler2012}.  In real glacier margins the ice can overhang, especially on steep bedrock features, violating our assumption of a single-valued surface elevation.  In fact, fractures, crevasses, and cliffs are common in ice margins, and even the bedrock can overhang.

However, all major numerical models ignore overhangs \cite{IsaacStadlerGhattas2015,Jouvetetal2008,LofgrenAhlkronaHelanow2022,WirbelJarosch2020}, and fractures do not occur within the viscous-fluid physics of Stokes models.  The current work follows these patterns.  Future models might allow fractures, supplementing momentum conservation with an advected damage variable and a stress-failure criterion \cite{PralongFunk2005}, and the emergent macroscopic margin shapes might suggest a regularity assumption and Sobolev space for surface elevations in the viscous theory.

We now attempt to choose $\cX$ based on mathematical considerations within the viscous theory.  Well-posedness of the fixed-domain Stokes problem \eqref{eq:glenstokes:weak}, and the surface trace bound in Lemma  \ref{lem:surfacetracebound}, generate a well-defined map from an admissible surface elevation $s\in \cK$ to the surface motion term $\bu|_s\cdot \bn_s$ in \eqref{eq:be:viearly}.  To clarify this map we use domain notation which emphasizes dependence on the elevation functions:
\begin{equation}
\Lambda(b,s) = \left\{(x,z)\,:\,b(x) < z < s(x)\right\} \subset \Omega \times \RR. \label{eq:domainfroms}
\end{equation}
This replaces the $t$-dependent notation in \eqref{eq:icydomain}; note $s(x)=s^n(x)\approx s(t_n,x)$ here.

For a given input surface $s$, the proposed map $s \mapsto \bu|_s\cdot \bn_s$ first solves the Stokes model \eqref{eq:glenstokes:weak} over $\Lambda(b,s)$, then evaluates the trace of $\bu$ along $\Gamma_s$, which is extended by zero \eqref{eq:defineus}, and finally it multiplies by $\bn_s$.

For this to be well-defined, $s\in\cK\subset\cX$ must be sufficiently-regular.  For $\cX$ we propose a Sobolev space $W^{1,\rr}(\Omega)$ because of the next Lemma.  Let $[L]>0$ be a representative \emph{horizontal} scale for $\Omega\subset\RR^2$, and for $\omega\in W^{1,\rr}(\Omega)$ define the norm
\begin{equation}
\|\omega\|_{W^{1,\rr}} = \left(\int_\Omega |\omega|^\rr + [L]^\rr \int_\Omega |\grad \omega|^\rr\right)^{1/\rr}. \label{eq:norm:Omega}
\end{equation}
(These are integrals over a fixed 2D domain, not over $\Lambda(b,s) \subset \RR^3$.)  We now show that if $\rr\ge 2$, $s\in W^{1,\rr}(\Omega)$, and $\Lambda(b,s)$ is Lipschitz, then the measurable function $\bu|_s\cdot \bn_s$ is in the dual space $(W^{1,\rr}(\Omega))'$.

\begin{lemma}[Surface motion is in the dual space] \label{lem:phibound}  Suppose $2 \le \rr \le \infty$ and that $s\in W^{1,\rr}(\Omega)$ satisfies $s\ge b$ and $s=b$ on $\partial\Omega$, and that $\Lambda(b,s)$ is a bounded Lipschitz domain.  Then there is a constant $C>0$, depending on $\|s\|_{W^{1,\rr}}$ and on the geometry of $\Lambda(b,s)$ through Lemma \ref{lem:surfacetracebound}, so that for all  $\omega\in W^{1,\rr}(\Omega)$,
\begin{equation}
\left|\int_\Omega \bu|_s\cdot \bn_s \,\omega\,dx\right| \le C\, \|\omega\|_{W^{1,\rr}}. \label{eq:phibound:early}
\end{equation}
\end{lemma}

\begin{proof}  We suppose $\rr<\infty$; the case $\rr=\infty$ is simpler because $|\bn_s| \le 1 + \|\grad s\|_{L^\infty}$ is then bounded.  Recall that $dS = |\bn_s|\,dx = \sqrt{1+|\grad s|^2}\,dx$ is the surface area element for $\Gamma_s \subset \partial \Lambda$.  Apply H\"older's inequalities:
\begin{align}
\Big|\int_\Omega \bu|_s\cdot \bn_s \, &\omega\,dx\Big| \le \int_\Omega \big|\bu|_s\big| |\bn_s|^{1/\pp} |\bn_s|^{1/\pp'} |\omega|\,dx \label{eq:phibound:zero} \\
    &\le \left(\int_\Omega \big|\bu|_s\big|^\pp |\bn_s|\,dx\right)^{1/\pp} \left(\int_\Omega |\bn_s| |\omega|^{\pp'} \,dx\right)^{1/\pp'} \notag \\
    &\le \left(\int_{\Gamma_s} |\bu|^\pp \,dS\right)^{1/\pp} \left(\int_\Omega |\bn_s|^\rr \,dx\right)^{1/(\pp'\rr)} \left(\int_\Omega |\omega|^{\pp'\rr'} \,dx\right)^{1/(\pp'\rr')}. \notag
\end{align}
From bound \eqref{eq:surfacetracebound} we may therefore write
\begin{equation}
\left|\int_\Omega \bu|_s\cdot \bn_s \, \omega\,dx\right| \le C \left(\int_\Omega \left(1+|\grad s|^2\right)^{\rr/2}\,dx\right)^{1/(\pp'\rr)} \|\omega\|_{L^{\pp'\rr'}}.
\end{equation}
Note that if $\alpha\ge 0$ then $(1+\alpha)^{\rr/2} \le 2^{(\rr-2)/2} (1+\alpha^{\rr/2})$, so, using generic constants,
\begin{align}
\left|\int_\Omega \bu|_s\cdot \bn_s \, \omega\,dx\right| &\le C \left(2^{(\rr-2)/2} \int_\Omega 1 + |\grad s|^\rr\,dx\right)^{1/(\pp'\rr)} \|\omega\|_{L^{\pp'\rr'}} \label{eq:phibound:one} \\
  &\le C \left(|\Omega| + [L]^{-\rr}\|s\|_{W^{1,\rr}}^\rr\right)^{1/(\pp'\rr)} \|\omega\|_{L^{\pp'\rr'}}. \notag
\end{align}
Since $2 \le \pp' \le \pp'\rr' < \infty$, by Sobolev's inequality \cite[Theorem 8.8, for example]{LiebLoss1997}, we have $\|\omega\|_{L^{\pp'\rr'}} \le C \|\omega\|_{W^{1,\rr}}$, thus \eqref{eq:phibound:early}.
\end{proof}

Because $\Omega\subset \RR^2$, if additionally $\rr>2$ then $W^{1,\rr}(\Omega) \hookrightarrow C(\bar\Omega)$.  Note that in Section \ref{sec:application} the FE approximation $s_h\approx s$ will be continuous and piecewise-linear.  Thus for $2 < \rr \le \infty$ we define from now on
\begin{equation}
\boxed{\cX = W^{1,\rr}(\Omega).} \label{eq:defineX}
\end{equation}
The admissible subset $\cK \subset \cX$ for surface elevations was already defined by \eqref{eq:be:admissible}.  By Lemma \ref{lem:surfacetracebound} the surface velocity $\bu|_s$ is well-defined as a measurable function, and its $L^\pp(\Gamma_s)$ norm is bounded.  By Lemma \ref{lem:phibound} the surface motion term $\bu|_s\cdot\bn_s$ is in $\cX'$.  Thus we officially define the nonlinear \emph{surface motion map} $\Phi:\cK \to \cX'$ by
\begin{equation}
\boxed{\Phi(s)[\omega] = -\int_\Omega \bu|_s\cdot\bn_s\,\omega\,dx.} \label{eq:definePhi}
\end{equation}

However, for well-posedness of VI \eqref{eq:be:viearly} we must conjecture that a norm of $\bu|_s$ is Lipschitz as a function of $s$.  Proving this would apparently require additional, detailed knowledge about Stokes problem \eqref{eq:glenstokes:weak}.  The dependence of scalar Laplacian-type equations on the domain is addressed in some literature, for example \cite{Daners2003}, but such results for Stokes problems are not known to this author.

\begin{conjecture}[Surface velocity is Lipschitz in the surface elevation] \label{conj:lipschitz}
There exists $2 < \rr \le \infty$ with the following two properties.  \emph{(i)} For every $s\in\cK$, Stokes problem \eqref{eq:glenstokes:weak} over $\Lambda(b,s)$ has a unique solution $(\bu,p)\in\mathcal{M}$, and bound \eqref{eq:surfacetracebound} holds, so the surface velocity $\bu|_s$ is well-defined over $\Omega$ by \eqref{eq:defineus}.  \emph{(ii)} If $R>0$ then there is $C(R)>0$ so that for $\sigma,s\in B_R \cap \cK = \{w \in \cK \,:\, \|w\|_{\cX} \le R\}$ we have
\begin{equation}
\big\|\bu|_\sigma - \bu|_s\big\|_{L^{\rr'}} \le C(R) \|\sigma-s\|_{\cX}. \label{eq:ulipschitz}
\end{equation}
\end{conjecture}

\begin{definition} \label{def:lipschitz}
For $R>0$ recall $B_R = \{v\in \cX\,:\,\|v\|\le R\}$.  We say that a map $f:\cK \to \cX'$ is \emph{Lipschitz on bounded subsets of $\cK$} if for every $R>0$ there is $C(R)>0$ so that
\begin{equation}
\|f(v)-f(w)\|_{\cX'} \le C(R) \|v-w\|_{\cX} \quad \text{ for all } v,w \in B_R \cap \cK.  \label{eq:liponbounded}
\end{equation}
\end{definition}

\begin{lemma} \label{lem:philipschitz}  If Conjecture \ref{conj:lipschitz} holds then $\Phi$ is Lipschitz on bounded subsets of $\cK$.
\end{lemma}

\begin{proof}  Assume without loss of generality that $R\ge \|b\|_\cX$, so $b\in B_R \cap \cK$.  Suppose $\sigma,s\in B_R \cap \cK$.  For any $\omega\in\cX$, in the following we add and subtract $\bu|_s \cdot \bn_\sigma$, and use the triangle inequality $|\bn_\sigma|=\left(1+|\grad \sigma|^2\right)^{1/2} \le 1 + |\grad \sigma|$:
\begin{align}
\Big|\Phi(\sigma)[\omega] - \Phi(s)[\omega]\Big| &\le \int_\Omega \Big|\bu|_\sigma - \bu|_s\Big| |\bn_\sigma| |\omega|\,dx + \int_\Omega \big|\bu|_s\big| |\bn_\sigma-\bn_s| |\omega|\,dx \\
    &\le \int_\Omega \Big|\bu|_\sigma - \bu|_s\Big| |\omega|\,dx + \int_\Omega \Big|\bu|_\sigma - \bu|_s\Big| |\grad \sigma| |\omega|\,dx \notag \\
    &\qquad\qquad + \int_\Omega \big|\bu|_s\big| |\grad \sigma-\grad s| |\omega|\,dx. \notag
\end{align}
By applying H\"older's inequality to each of these integrals we have
\begin{align}
\int_\Omega \Big|\bu|_\sigma - \bu|_s\Big| |\omega|\,dx &\le \big\|\bu|_\sigma - \bu|_s\big\|_{L^{\rr'}} \|\omega\|_{L^{\rr}} \le \|\bu|_\sigma - \bu|_s\|_{L^{\rr'}} \|\omega\|_{\cX}, \label{eq:philipschitz:1} \\
\int_\Omega \Big|\bu|_\sigma - \bu|_s\Big| |\grad \sigma| |\omega|\,dx &\le \left(\int_\Omega \Big|\bu|_\sigma - \bu|_s\Big|^{\rr'} |\omega|^{\rr'}\, dx\right)^{1/\rr'} \|\grad \sigma\|_{L^\rr} \label{eq:philipschitz:2} \\
    &\le [L]^{-1} \big\|\bu|_\sigma - \bu|_s\big\|_{L^{\rr'}} \|\sigma\|_{\cX} \|\omega\|_{L^\infty}, \notag \\
\int_\Omega \big|\bu|_s\big| |\grad \sigma-\grad s| |\omega|\,dx &\le \left(\int_\Omega \big|\bu|_s\big|^{\rr'} |\omega|^{\rr'}\, dx\right)^{1/\rr'} \|\grad \sigma- \grad s\|_{L^\rr}  \label{eq:philipschitz:3} \\
    &\le [L]^{-1} \big\|\bu|_s - \bzero\big\|_{L^{\rr'}} \|\sigma-s\|_{\cX} \|\omega\|_{L^\infty}. \notag
\end{align}
Note that $\bu|_b=\bzero$; there is no flow when there is no glacier.  Because $\rr>2$, Sobolev's inequality gives $\|\omega\|_{L^\infty} \le c_\infty \|\omega\|_\cX$ for some $c_\infty>0$.  Now apply Conjecture \ref{conj:lipschitz} to each of \eqref{eq:philipschitz:1}--\eqref{eq:philipschitz:3}:
\begin{equation}
\big|\Phi(\sigma)[\omega] - \Phi(s)[\omega]\big| \le \tilde C(R) \left(1 + c_\infty [L]^{-1} \left(\|\sigma\|_{\cX} + \|s - b\|_{\cX}\right)\right) \|\sigma-s\|_{\cX} \|\omega\|_{\cX}.
\end{equation}
Recalling that $b\in B_R\cap \cK$, use the triangle inequality again to conclude with $C(R) = \tilde C(R) \left(1 + 3 c_\infty [L]^{-1} R\right)$.
\end{proof}

At this point we have stated the tools needed to make the continuum free-boundary model for a backward-Euler time-step \eqref{eq:be:ske} mathematically precise.

\begin{definition} For $\Delta t>0$, $s\in\cK$, and $\omega\in\cX$ define the nonlinear (Stokes) \emph{geometry update operator} $F_{\Delta t}:\cK\to\cX'$ as
\begin{equation}
\boxed{F_{\Delta t}(s)[\omega] = \Delta t\,\Phi(s)[\omega] + \int_\Omega s \omega.}  \label{eq:be:Fdefine}
\end{equation}
Assume that $\ell^n$ is defined as an integrable function from the data by \eqref{eq:be:source}, and by slight abuse of notation write $\ell^n[\omega] = \int_\Omega \ell^n \omega$ so that $\ell^n\in \cX'$.  We say that the surface elevation $s=s^n \in \cK$ solves the (weak-form) \emph{implicit time-step problem} if
\begin{equation}
\boxed{F_{\Delta t}(s)[\sigma-s] \ge \ell^n[\sigma-s] \quad \text{for all } \sigma \in \cK.} \label{eq:be:vi}
\end{equation}
\end{definition}

By Lemma \ref{lem:phibound} and Conjecture \ref{conj:lipschitz}(i), $\Phi$ and $F_{\Delta t}$ are well-defined, and if Conjecture \ref{conj:lipschitz} holds then they are Lipschitz on bounded subsets.  The reader should confirm that the weak-form \eqref{eq:be:vi} then merely rewrites \eqref{eq:be:viearly}.  This VI gives a precise meaning to strong-form NCP \eqref{eq:be:ncp}.

The boxed definitions and inequalities above, namely \eqref{eq:be:admissible}, \eqref{eq:be:source}, \eqref{eq:defineX}, \eqref{eq:definePhi}, \eqref{eq:be:Fdefine}, and \eqref{eq:be:vi}, form a continuum model for a single time-step of glacier evolution, based on the solution of Stokes problem \eqref{eq:glenstokes:weak} over the domain $\Lambda=\Lambda(b,s)$ defined in \eqref{eq:domainfroms}.  Although Conjecture \ref{conj:lipschitz} is used in the construction of the operator $F_{\Delta t}$, the whole model is precise enough to be mathematically analyzed.

A future existence result for \eqref{eq:be:vi} is possible, but we will regularize the model in the next section.  Based on numerical evidence, well-posedness for this regularization is credible.  Sections \ref{sec:abstractestimate} and \ref{sec:application} then prove an error estimation theorem for an FE approximation of the regularized model.

\section{Conjectural well-posedness for the regularized problem} \label{sec:conjectural}

Theorem \ref{thm:glacierapp}, the main result to come, bounds numerical error.  It bounds the FE surface elevation error by applying Theorem \ref{thm:abstractestimate} (next section), which uses the concepts of monotonicity \cite{Minty1963} and coercivity \cite[Chapter III]{KinderlehrerStampacchia1980}.

\begin{definition} \label{def:monotonecoercive}
Suppose $\cX$ is a Banach space, with norm $\|\cdot\|$, and $\cK\subset \cX$ is closed and convex.  A generally-nonlinear operator $f:\cK \to \cX'$ is called \emph{monotone} if
\begin{equation}
\left(f(v)-f(w)\right)[v-w] \ge 0 \qquad \text{for all } v,w \in \cK, \label{eq:monotone}
\end{equation}
and \emph{strictly monotone} if equality implies $v=w$.  It is \emph{coercive} if there is $w\in \cK$ so that $\left(f(v)-f(w)\right)[v-w]/\|v-w\| \to +\infty$ for $v \in \cK$ as $\|v\| \to +\infty$.  It is \emph{$\qq$-coercive} for $\qq>1$ if there exists $\alpha>0$ such that
\begin{equation}
\left(f(v)-f(w)\right)[v-w] \ge \alpha \|v-w\|^\qq \qquad \text{for all } v,w \in \cK. \label{eq:qcoercive}
\end{equation}
\end{definition}

Note that $\qq$-coercivity \cite{Bueler2021conservation} implies coercivity and strict monotonicity, and that $2$-coercivity is sometimes called \emph{strong monotonicity} \cite{Chow1989}.

Theorem \ref{thm:abstractestimate} applies to operators $f:\cK\to\cX'$ which are $\qq$-coercive.  We would like to apply it to $F_{\Delta t}$, defined by \eqref{eq:be:Fdefine}, over the admissible subset $\cK \subset W^{1,\rr}(\Omega)$.  However, examples\footnote{These apparently-new examples also negatively resolve a previously open question of uniqueness, for elevation-independent SMB, in the SIA.} in Appendix \ref{app:noncoercive} show that, for any bedrock $b$ possessing strict local minima in the interior of $\Omega$, the surface motion map $\Phi$ \eqref{eq:definePhi} alone is not $\qq$-coercive, nor strictly monotone, over $W^{1,\rr}(\Omega)$ for any $\rr$.  (However, since $F_{\Delta t}$ adds the integral $\int_\Omega s\omega$ to the surface motion term $\Delta t\,\Phi(s)[\omega]$, lack of coercivity for $\Phi$ does not dis-prove it for $F_{\Delta t}$.)

The numerical experiments in Appendix \ref{app:numerical} show that the surface motion operator $\Phi$ is nonetheless \emph{nearly} $\qq$-coercive.  In these experiments thousands of high-resolution surface pairs $\sigma,s\in \cK\subset W^{1,\rr}(\Omega)$ were generated from glacier simulations over $\Omega \subset \RR^1$.  Coercivity ratios for $\Phi$ were computed with exponent $\rr=\qq=4$,
\begin{equation}
\frac{(\Phi(\sigma) - \Phi(s))[\sigma - s]}{\|\sigma-s\|_{W^{1,4}}^{4}}, \label{eq:ratios}
\end{equation}
using surface velocities from a numerical Stokes model when evaluating $\Phi$.  By definition \eqref{eq:qcoercive}, if $\Phi$ were $4$-coercive over $\cX=W^{1,4}(\Omega)$, and if computations were exact, then these ratios would all exceed a constant $\alpha>0$.  Appendix \ref{app:numerical} shows that pairs yielding negative ratios do exist, about $10\%$ of samples, but that their magnitudes are an order of magnitude smaller than those of the positive ratios.  Standard presentations of Stokes models for glaciers give no direct reason to expect these ratios should prefer positive values, but we see that they decidedly do.

To suggest why coercivity might occur, we can rewrite the surface motion $\Phi(s)$ in a different manner.  Extend $\omega\in \cX$ to $\Lambda(b,s)$ as $\tilde\omega(x,z)=\omega(x)$, i.e.~independent of $z$, and define the glaciological ice flux $\bq(s) = \int_b^s \bU\,dz$ \cite{SchoofHewitt2013}, where $\bU=(u_1,u_2)$ is the horizontal velocity in $\bu=(u_1,u_2,w)$.  Definition \eqref{eq:definePhi}, the non-sliding assumption ($\bu|_{\Gamma_b}=\bzero$), the divergence theorem, and incompressibility yield
\begin{align}
\Phi(s)[\omega] &= -\int_{\Omega} \bu|_s \cdot \bn_s\,\omega\,dx = -\int_{\partial\Lambda(b,s)} \tilde\omega \bu\cdot {\hat\bn}\,dS = - \int_{\Lambda(b,s)} \Div\left(\tilde\omega \bu\right)  \label{eq:fluxformPhi} \\
  &= - \int_{\Lambda(b,s)} \bu\cdot \grad \tilde\omega = -\int_\Omega \bq(s) \cdot \grad \omega, \notag
\end{align}
where $\hat\bn$ is the outward unit normal and $dS$ is the area element along $\partial\Lambda(b,s)$.  Now, the natural glaciological belief is that ice flows downhill, which is the informal association $\bq(s) \sim - \grad s$, in which case the last expression in \eqref{eq:fluxformPhi} implies $(\Phi(\sigma) - \Phi(s))[\sigma - s] \sim \int_\Omega |\grad s - \grad \sigma|^2$, which is $2$-coercivity in the $H^1$ norm.  Observe that the rheology of ice, Glen's law \eqref{eq:glen}, does not enter into this informal reasoning, so it cannot be the full story regarding coercivity.

\newcommand{\bqSIA}{\bq^{\text{SIA}}}
\newcommand{\PhiSIA}{\Phi^{\text{SIA}}}

To enhance coercivity of the surface motion $\Phi$, coming from the Stokes surface velocity, consider the ice flux in the isothermal SIA \cite{GreveBlatter2009}:
\begin{equation} \label{eq:qsia}
\bqSIA(s) = - \Gamma (s-b)^{\nn+2} |\grad s|^{\nn-1} \grad s.
\end{equation}
Here $\nn$ is Glen's exponent and $\Gamma>0$ is an ice softness constant computable from the viscosity scale $\nu_\pp$ in \eqref{eq:glen}.  Formula \eqref{eq:qsia} does \emph{not} describe the ice flux in the full Stokes model; we will only use it for regularization.  Because it degenerates at the glacier margin as the ice thickness goes to zero, we replace $s-b$ with a fixed value $\eta_0>0$.  Via an argument like \eqref{eq:fluxformPhi}, the corresponding surface motion has the form
\begin{equation}
\PhiSIA(s)[\omega] = \Gamma \eta_0^{\nn+2} \int_\Omega |\grad s|^{\nn-1} \grad s \cdot \grad \omega \label{eq:phisia}
\end{equation}
for $\omega \in \cX$.  Now for $0<\eps \ll 1$, the regularization is this convex combination:
\begin{equation} \label{eq:regularizedPhi}
\Phi^\eps(s)[\omega] = (1-\eps) \Phi(s)[\omega] + \eps \PhiSIA(s)[\omega].
\end{equation}
This regularizes the Stokes surface motion \eqref{eq:definePhi} using a small amount from formula \eqref{eq:phisia}.  If $\Phi$ is Lipschitz on bounded subsets of $\cK$ (Lemma \ref{lem:philipschitz}) then so is $\Phi^\eps$, at least when $\rr\ge\nn+1$.  This follows from the inequality $\big||\xi|^{\nn-1}\xi - |\zeta|^{\nn-1}\zeta\big| \le \nn (|\xi|+|\zeta|)^{\nn-1}|\xi-\zeta|$ for $\xi,\zeta\in\RR^2$, and H\"older's inequality with exponents $\rr/(\nn-1)$, $\rr$, $\rr$, whose reciprocals sum to $(\nn+1)/\rr \le 1$.  It is reasonable to suppose $\Phi^\eps$ is a $4$-coercive operator over $\cK\subset W^{1,4}(\Omega)$ because of properties of the $4$-Laplacian operator \cite[for example]{ChoeLewis1991}; with $\nn=3$ and $\eta_0>0$, $\PhiSIA$ in \eqref{eq:phisia} is a multiple of the $4$-Laplacian, which is $4$-coercive over $\cK$ because differences $\sigma-s$ of elements of $\cK$ have zero trace on $\partial\Omega$.

The numerical experiments in Appendix \ref{app:numerical} use values $\eps=0.01$ and $\eta_0=500$ meters in \eqref{eq:phisia} and \eqref{eq:regularizedPhi}.  The observed effect is that negative values of the regularized coercivity ratios, i.e.~\eqref{eq:ratios} with $\Phi^\eps$ in place of $\Phi$, disappear at high spatial resolution, and the distribution of ratios is bounded away from zero.  However, we are not able to \emph{prove} that $\Phi^\eps$ in \eqref{eq:regularizedPhi} is $4$-coercive.  A proof would seem to depend on novel insights into the glaciological Stokes problem.  Thus the numerical analysis in Section \ref{sec:application} will be based on the following second conjecture.

\begin{conjecture}[Regularized surface motion $\Phi^\eps(s)$ is $4$-coercive on $\cK\subset W^{1,4}$] \label{conj:regcoercive}
Let $\cX = W^{1,4}(\Omega)$ and $\nn=3$ in \eqref{eq:phisia}.  Fix $b\in C^1(\bar\Omega)$ and let $\cK=\{\sigma\in\cX\,:\,\sigma\ge b \text{ and } \sigma|_{\partial\Omega}=b|_{\partial\Omega}\}$.  Suppose Conjecture \ref{conj:lipschitz} holds for $\rr=4$.  There exists $\eps \in (0,1)$ and $\eta_0>0$, and $\alpha>0$, so that
\begin{equation}
\left(\Phi^\eps(\sigma) - \Phi^\eps(s)\right)[\sigma-s] \ge \alpha \|\sigma-s\|_{\cX}^4 \qquad \text{for all } \sigma,s\in\cK. \label{eq:regcoercive}
\end{equation}
\end{conjecture}

This Conjecture suffices for well-posedness of the regularized VI problem.

\begin{theorem} \label{thm:regularizedwellposed}  Assume Conjectures \ref{conj:lipschitz} and \ref{conj:regcoercive}.  Suppose that $s^{n-1}\in\cK$ and define the source term $\ell^n \in \cX'$ as in \eqref{eq:be:source}.  For $\Delta t>0$, and $\eps,\eta_0$ as in Conjecture \ref{conj:regcoercive}, define
\begin{equation}
F^\eps_{\Delta t}(s)[\omega] = \Delta t\,\Phi^\eps(s)[\omega] + \int_\Omega s \omega, \label{eq:regularizedF}
\end{equation}
the regularized update operator.  This operator is Lipschitz on bounded subsets and $4$-coercive, thus there exists a unique surface elevation $s\in\cK$ satisfying the VI problem%
\begin{equation}
F^\eps_{\Delta t}(s)[\sigma-s] \ge \ell^n[\sigma-s] \quad \text{for all } \sigma \in \cK. \label{eq:regularizedvi}
\end{equation}
\end{theorem}

\begin{proof}  The Lipschitz property follows from Lemma \ref{lem:philipschitz} and the remark after \eqref{eq:regularizedPhi}, and $4$-coercivity follows from Conjecture \ref{conj:regcoercive}.  In particular, $F^\eps_{\Delta t}$ is $4$-coercive over $\cK$ with constant at least $\alpha\Delta t>0$, and thus it is also coercive and strictly-monotone (Definition \ref{def:monotonecoercive}).  Corollary III.1.8 of \cite{KinderlehrerStampacchia1980} shows existence for \eqref{eq:regularizedvi}, and strict monotonicity gives uniqueness.
\end{proof}

Theorem \ref{thm:regularizedwellposed} addresses only the well-posedness of a single time-step, over $[t_{n-1},t_n]$.  Its conclusion is not sufficient to show well-posedness of the time-dependent VI problem corresponding to NCP \eqref{eq:ncp}, nor to show that implicit steps converge in the $\Delta t\to 0$ limit.  However, it is a first mathematical step in these directions.  Time-stepping numerical models for the evolution of glacier geometry, using Stokes dynamics, are already in common use.  Practitioners apparently expect that each time-step problem, implicit or not, is well-posed, and that computed surface elevations will converge to continuum solutions in some well-behaved manner.

\section{Abstract error estimate for finite element approximation of VIs} \label{sec:abstractestimate}

In this Section we consider the FE approximation of an abstract VI problem in a Banach space.  We will return to glaciological problem \eqref{eq:regularizedvi} in Section \ref{sec:application}.

Let $\cX$ be a real, reflexive Banach space with norm $\|\cdot\|$ and topological dual space $\cX'$.  Denote the dual pairing of $\ell \in \cX'$ and $v\in\cX$ by $\ell[v]$, and define $\|\ell\|_{\cX'} = \sup_{\|v\|=1} \big|\ell[v]\big|$.  Let $\cK \subset \cX$ be a nonempty, closed, and convex subset, the admissible set.  For a continuous and (generally) nonlinear operator $f:\cK \to \cX'$, and a source functional $\ell\in \cX'$, the VI problem is to find $u\in \cK$ such that
\begin{equation}
f(u)[v-u] \ge \ell[v-u] \quad \text{for all } v\in \cK. \label{eq:vi}
\end{equation}
Problem \eqref{eq:regularizedvi} is in this form, while the best known example for \eqref{eq:vi} is the obstacle problem for the Laplacian operator; see \cite{Ciarlet2002,Evans2010,KinderlehrerStampacchia1980} for theory and FE analysis.

If $f:\cK \to \cX'$ is monotone and coercive (Definition \ref{def:monotonecoercive}), and continuous on finite-dimensional subspaces, then problem \eqref{eq:vi} has a solution \cite[Corollary III.1.8]{KinderlehrerStampacchia1980}, and if $f$ is strictly monotone then the solution is unique.  If $f$ is $\qq$-coercive then it is coercive and strictly monotone, and if $f$ is Lipschitz on bounded subsets (Definition \ref{def:lipschitz}) then it is continuous.  Recall that $f(u)-\ell \in \cX'$ is generally nonzero if $u$ solves \eqref{eq:vi}.  Under sufficient regularity assumptions, an NCP strong form follows from \eqref{eq:vi}.

Now suppose $\cX_h \subset \cX$ is a finite-dimensional subspace, an FE space of continuous, piecewise-polynomial functions defined on a mesh.  Let $\cK_h\subset \cX_h$ be a closed and convex subset, and suppose $f_h:\cK_h\to\cX'$ is an approximation of $f$.  Generally $\cK_h \nsubset \cK$---we will \emph{not} assume initially that FE states are continuum admissible---and generally $f_h\ne f$ because of quadrature and other approximations.  In fact, in standard FE Stokes-type glacier simulations both $\cK_h \nsubset \cK$ and $f_h\ne f$ apply for practical reasons, though in Section \ref{sec:application} we arrange that $\cK_h\subset\cK$ by construction of the meshed bed elevation.

The FE method for \eqref{eq:vi} is the finite-dimensional VI problem
\begin{equation}
f_h(u_h)[v_h-u_h] \ge \ell[v_h-u_h] \quad \text{for all } v_h\in \cK_h. \label{eq:fe:vi}
\end{equation}
We will assume that this problem has a solution $u_h\in\cK_h$.  The following estimation theorem extends the well-known result of Falk \cite{Falk1974}; see also Theorem 5.1.1 in \cite{Ciarlet2002}.  We assume that $f$ is $\qq$-coercive and Lipschitz on some superset of $\cK$ which also contains the FE solution $u_h$; this technical assumption permits a general result.  However, we do \emph{not} assume any of the following: $\cK_h \subset \cK$, $f$ is linear, $f_h=f$, $f_h$ is continuous, $f_h$ is $\qq$-coercive, or even that the solution of \eqref{eq:fe:vi} is unique.

\begin{theorem} \label{thm:abstractestimate}  Suppose $u\in\cK$ solves \eqref{eq:vi} and $u_h\in\cK_h$ solves \eqref{eq:fe:vi}.  For $\qq>1$, with conjugate exponent $\qq'=\qq/(\qq-1)$, assume that $f$ is $\qq$-coercive, with constant $\alpha>0$, on a set which contains $\cK$ and the numerical solution $u_h$.  Assume $f$ is Lipschitz on bounded subsets of the same set, and let $L(R)$ denote the Lipschitz constant in \eqref{eq:liponbounded} on $B_R$.  Let $R_h=\max\{\|u\|,\|u_h\|\}$.  Then
\begin{align}
\frac{\alpha}{2} \|u-u_h\|^\qq &\le \left(f(u)-\ell\right)[u-u_h] \label{eq:abstractestimate} \\
   &\quad\, + \inf_{v_h\in\cK_h} \bigg\{ \left(f(u)-\ell\right)[v_h-u] + \delta_h[u_h-v_h] + C \|v_h - u\|^{\qq'} \bigg\}, \notag
\end{align}
where $\delta_h = f(u_h)-f_h(u_h) \in \cX'$ is a consistency error at the discrete solution, and
\begin{equation}
C = \frac{1}{\qq'} \left(\frac{2}{\qq\,\alpha}\right)^{1/(\qq-1)} L(R_h)^{\qq'}. \label{eq:abstractconstant}
\end{equation}
\end{theorem}

\begin{proof}  It follows from $\qq$-coercivity of $f$ that
\begin{align}
\alpha \|u-u_h\|^\qq &\le \left(f(u)-f(u_h)\right)[u-u_h] \notag \\
  &= f(u)[u] + f(u_h)[u_h] - f(u)[u_h] - f(u_h)[u] \notag \\
  &= f(u)[u] + f_h(u_h)[u_h] - f(u)[u_h] - f(u_h)[u] + \left(f(u_h)-f_h(u_h)\right)[u_h]. \label{eq:abstract:one}
\end{align}
Next, for arbitrary $v\in\cK$ and $v_h\in\cK_h$, rewrite \eqref{eq:vi} as $f(u)[u]     \le f(u)[v] + \ell[u-v]$ and \eqref{eq:fe:vi} as $f_h(u_h)[u_h] \le f_h(u_h)[v_h] + \ell[u_h-v_h]$.  Then
\begin{align}
\alpha \|u-u_h\|^\qq &\le f(u)[v] + \ell[u-v] + f_h(u_h)[v_h] + \ell[u_h-v_h] \label{eq:abstract:two} \\
  &\qquad - f(u)[u_h] - f(u_h)[u] + \left(f(u_h)-f_h(u_h)\right)[u_h] \notag \\
  &= \left(f(u)-\ell\right)[v-u_h] + \left(f(u)-\ell\right)[v_h-u] \notag \\
  &\qquad + \left(f(u)-f(u_h)\right)[u-v_h] + \left(f(u_h)-f_h(u_h)\right)[u_h - v_h]. \notag
\end{align}
Since $u,u_h\in B_{R_h}$, by the Lipschitz assumption
\begin{equation}
\left(f(u)-f(u_h)\right)[u-v_h] \le L(R_h) \|u-u_h\|\|u-v_h\|. \label{eq:abstract:three}
\end{equation}
Noting $1<\qq<\infty$, use Young's inequality with $\eps>0$ \cite[Appendix B.2]{Evans2010} on \eqref{eq:abstract:three}:
\begin{align}
\alpha \|u-u_h\|^\qq &\le \left(f(u)-\ell\right)[v-u_h] + \left(f(u)-\ell\right)[v_h-u]  \label{eq:abstract:four} \\
  &\qquad + L(R_h) \left(\eps\|u-u_h\|^\qq + \tilde C(\eps) \|u-v_h\|^{\qq'}\right) \notag \\
  &\qquad + \left(f(u_h)-f_h(u_h)\right)[u_h - v_h], \notag
\end{align}
where $\tilde C(\eps) = (\eps \qq)^{-\qq'/\qq} {\qq'}^{-1}$.  Choose $\eps = \alpha/(2L(R_h))$ and subtract:
\begin{align}
\frac{\alpha}{2} \|u-u_h\|^\qq &\le \left(f(u)-\ell\right)[v-u_h] + \left(f(u)-\ell\right)[v_h-u]  \label{eq:abstract:five} \\
  &\qquad + L(R_h) \tilde C(\eps) \|u-v_h\|^{\qq'} + \left(f(u_h)-f_h(u_h)\right)[u_h - v_h]. \notag
\end{align}
For the first term note that $(f(u)-\ell)[v-u_h] = (f(u)-\ell)[v-u] + (f(u)-\ell)[u-u_h] \ge (f(u)-\ell)[u-u_h]$ by \eqref{eq:vi}; the VI implies that the infimum over $v\in\cK$ is achieved by $u$.  Take infimum over $\cK_h$ to show \eqref{eq:abstractestimate}, noting that $L(R_h)\tilde C(\alpha/(2L(R_h)))$ equals $C$ in \eqref{eq:abstractconstant} because $1+\qq'/\qq=\qq'$.
\end{proof}

In a unilateral obstacle problem over scalar functions $\cX=W^{1,\rr}(\Omega)$, we may rewrite portions of bound \eqref{eq:abstractestimate} using integration.  In this case the residual $f(u)-\ell\in \cX'$ is represented by a positive Borel measure $\mu_u$ \cite[Theorem 6.22]{LiebLoss1997} which, by the proof of Theorem II.6.9 in \cite{KinderlehrerStampacchia1980}, is supported in the active set $A_u=\{x\in\Omega\,:\,u(x)=b(x)\}$.  Then $(f(u) - \ell)[w] = \int_{A_u} w\,d\mu_u$ for $w\in \cX$ with zero trace.

If $\cK_h \subset \cK$ then the first term in \eqref{eq:abstractestimate} is nonpositive and it can be replaced by zero; substitute $v=u_h$ into \eqref{eq:vi} to get $(f(u)-\ell)[u-u_h] \le 0$.

However, when $\cK_h \not\subset \cK$ the first term can be positive.  This can happen even in an obstacle problem with an arbitrarily-smooth obstacle.  For $\cK$ defined by \eqref{eq:be:admissible}, assume that $\Omega$ is polygonal, so that a mesh covers it exactly, and suppose $b_h=\pi_h b$ is the FE interpolant of $b$, and define $\cK_h=\{v_h \in \cX_h\,:\,v_h\ge b_h \text{ and } v_h|_{\partial\Omega}=b_h|_{\partial\Omega}\}$.  While $b_h(x_j)=b(x_j)$ for the interpolation nodes $x_j$, generally $b_h(x) \ge b(x)$ will not hold for all $x\in\Omega$ (Figure \ref{fig:nonadmissible}; compare \cite[Figure 5.1.3]{Ciarlet2002}); nodal admissibility does not imply admissibility.  Now if $b>b_h$ in the active set $A_u$ then generally $u>u_h$ also, and thus $(f(u)-\ell)[u-u_h] = \int_{A_u} (u-u_h)\,d\mu_u >0$.

\begin{figure}[ht]
\begin{center}
\begin{tikzpicture}[scale=1.1, domain=0.0:4.0, samples=200]
  \draw[black,thin,->] (-0.3,0.0) -- (4.5,0.0) node [xshift=2mm] {$x$};
  \draw plot (\x, {1.5 + 0.8 * cos(1.8*\x r)});
  \node[yshift=-2mm] at (2.0, 0.7) {$b$};
  \newcommand{\xlist}{0.0, 1.0, 1.4, 2.4, 3.0, 4.0}
  \foreach \x [remember=\x as \lastx] in \xlist {
      \draw (\x, 0.0) circle (1.5pt);
      \filldraw (\x, {1.5 + 0.8 * cos(1.8*\x r)}) circle (1.5pt);
      \draw[dashed] (\lastx, {1.5 + 0.8 * cos(1.8*\lastx r)}) -- (\x, {1.5 + 0.8 * cos(1.8*\x r)});
  }
  \node[yshift=3mm] at (3.0, 0.0) {$x_i$};
  \node[xshift=3mm, yshift=-7mm] at (0.0, 2.3) {$b_h$};
\end{tikzpicture} \quad \begin{tikzpicture}[scale=1.1, domain=0.0:4.0, samples=200]
  \draw[black,thin,->] (-0.3,0.0) -- (4.5,0.0) node [xshift=2mm] {$x$};
  \draw plot (\x, {1.5 + 0.8 * cos(1.8*\x r)});
  \newcommand{\xlist}{0.0, 1.0, 1.4, 2.4, 3.0, 4.0}
  \newcommand{\xymaxlist}{0.0/2.3, 1.0/2.3, 1.4/1.3182, 2.4/2.0078, 3.0/2.3, 4.0/2.3}
  \foreach \x/\y in \xymaxlist {
      \draw (\x, 0.0) circle (1.5pt);
      \filldraw (\x, \y) circle (1.5pt);
  }
  \node[yshift=3mm] at (3.0, 0.0) {$x_i$};
  \draw[dashed] (0.0,2.3) -- (1.0,2.3) -- (1.4,1.3182) -- (2.4,2.0078) -- (3.0,2.3) -- (4.0,2.3);
  \node[yshift=-2mm] at (2.0, 0.7) {$b$};
  \node[xshift=5mm, yshift=-3mm] at (1.0, 2.3) {$b_h$};
\end{tikzpicture}
\end{center}
\caption{Left: If $b_h=\pi_h b$ is the interpolant of $b$ then generally $\cK_h\nsubset\cK$.  Right: Defining $b_h=R_h^{\oplus} b$, using \eqref{eq:monotoneop}, implies $b_h\ge b$.}
\label{fig:nonadmissible}
\end{figure}

For unilateral obstacle problems and $P_1$ or $Q_1$ elements we may force $b_h \ge b$ by using a monotone restriction operator \cite{BuelerFarrell2024}, defined as follows.  For a node $x_i$ let $N_i$ be the union of elements adjacent to $x_i$.  Define $R_h^{\oplus}b\in\cX_h$ by
\begin{equation}
(R_h^{\oplus}b)(x_i) = \sup_{x \in N_i} b(x). \label{eq:monotoneop}
\end{equation}
If $b_h=R_h^{\oplus} b$ then $b_h(x)\ge b(x)$ for all $x\in \Omega$ (Figure \ref{fig:nonadmissible}; right).

Generally $b_h\ne b$ along $\partial\Omega$, which blocks $\cK_h\subset\cK$, so a boundary-modified version of \eqref{eq:monotoneop} sets $b_h(x_i)=b(x_i)$ at boundary nodes.  Then $\cK_h\subset\cK$ holds if, for example, $b$ is affine on each element meeting $\partial\Omega$.

The appearance of a consistency error $\delta_h = f(u_h)-f_h(u_h)$ in bound \eqref{eq:abstractestimate}, evaluated at $u_h\in \cK_h$, suggests Strang's first lemma \cite[Theorem 4.1.1]{Ciarlet2002}, but here in a nonlinear case.  In the glacier application this term will be related to the maximum error made in approximating the (Stokes) ice velocity; see Theorem \ref{thm:glacierapp}.

Falk's original VI bound \cite[Theorem 1]{Falk1974} is recovered from \eqref{eq:abstractestimate} in the Hilbert-space case with a bilinear, uniformly-elliptic operator and $f_h=f$, by bounding each dual pairing involving $f(u)-\ell$ with a norm of the residual, for example in $L^2(\Omega)$.  Doing so, however, loses the active-set information implicit in the residual $f(u)-\ell \in \cX'$.

For an obstacle problem, consider the case where the inactive set for $u$ is all of $\Omega$, thus $f(u)-\ell=0$.  (For glaciers this situation occurs when the entire domain $\Omega$ is covered by ice.)  If also $f=f_h$ then bound \eqref{eq:abstractestimate} reduces to its final term.  Thus we get a Cea-type bound $\|u-u_h\| \le C \inf_{v_h\in\cK_h} \|v_h-u\|^{1/(\qq-1)}$, in a Banach space and for a $\qq$-coercive nonlinear operator, which is Cea's lemma \cite[Theorem 2.4.1]{Ciarlet2002} when $\qq=2$.

In summary, Theorem \ref{thm:abstractestimate} generalizes Falk's bound, Strang's lemma, and Cea's lemma.  It is well-suited for the glacier application in the next section.

\section{Application to numerical glacier models} \label{sec:application}

Now we can apply the above theory to a time-step of an evolving-surface, Stokes-based glacier simulation.  Our approach essentially defines the phrase ``conforming FE method'' for such models.  We will combine the major threads so far: a surface velocity bound for the glaciological Stokes problem on a fixed domain (Section \ref{sec:stokes}), conjectural Lipschitz continuity and $4$-coercivity of the regularized surface-motion operator (Sections \ref{sec:model} and \ref{sec:conjectural}), and an abstract error estimate for FE solutions of unilateral obstacle problems (Section \ref{sec:abstractestimate}).

\begin{figure}[ht]
\begin{center}
\includegraphics[width=0.7\textwidth]{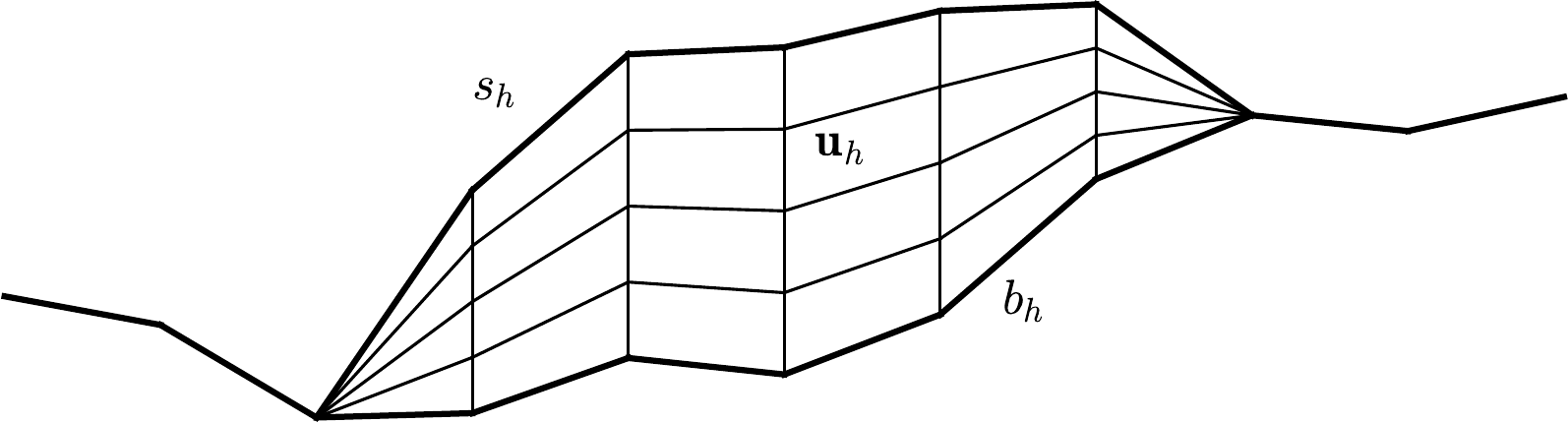}
\end{center}
\caption{The numerical Stokes problem computes the velocity $\bu_h$ on a 3D polyhedral domain between $P_1$ functions $b_h$ and $s_h$.  Evaluating $F_{\Delta t}^{\eps,h}$, formula \eqref{eq:fe:regularizedF}, uses the surface trace $\bu_h|_{s_h}$.}
\label{fig:fe:domain}
\end{figure}

Practical models for glaciers, based on unstructured meshes in the 2D map-plane, use polyhedral domain geometry for the 3D Stokes problem \cite{Brinkerhoff2023,IsaacStadlerGhattas2015,Jouvetetal2008,Lengetal2012}.  For this reason we assume that the discrete bed and surface elevations $b_h,s_h$ are from the continuous $P_1$ FE space $\cX_h$ over a triangulation $\cT_h$ of the domain $\Omega \subset \RR^2$ (Figure \ref{fig:fe:domain}).  Furthermore the 3D mesh might be extruded vertically from $\cT_h$, as shown in Figure \ref{fig:fe:domain}, but we do not specifically assume so.

We will assume that the original bed elevation $b$ is well-behaved, with $b\in C^1(\bar\Omega)$.  (In practice the data $b$ is a high-resolution map derived from ice-penetrating radar \cite{Morlighemetal2017}, typically over a much finer mesh than $\cT_h$; here it is assumed to be smooth.)  For the discrete bed $b_h\in\cX_h$ we suppose $b_h=b$ along $\partial\Omega$ and that $b_h\ge b$, for instance by a boundary-modified monotone restriction (Section \ref{sec:abstractestimate}), which suffices if $b$ is affine on each element meeting $\partial\Omega$.  The fact that $b_h$ is above the data $b$ eliminates a bound term when Theorem \ref{thm:abstractestimate} is applied below.  Following \eqref{eq:be:admissible}, we define $\cK_h=\{\sigma_h \in \cX_h\,:\,\sigma_h|_{\partial\Omega}=b_h|_{\partial\Omega} \text{ and } \sigma_h\ge b_h\}$, so that $\cK_h\subset\cK \subset \cX$ where $\cX = W^{1,4}(\Omega)$.

Along with these geometrical and discrete-space choices, recall the regularized VI problem \eqref{eq:regularizedvi}, with $F^\eps_{\Delta t}$ defined by \eqref{eq:regularizedPhi} and \eqref{eq:regularizedF}.  Numerical approximation of $F^\eps_{\Delta t}(s_h)$ will use the numerical solution of the power-law Stokes problem \eqref{eq:glenstokes:weak} over a meshed 3D domain $\Lambda(b_h,s_h)$---see definition \eqref{eq:domainfroms}---the lower and upper surfaces of which are given by admissible FE functions $b_h\in\cX_h$ and $s_h \in \cK_h$, respectively.  Solvability of this Stokes problem using inf-sup stable elements is addressed by \cite{Belenkietal2012,JouvetRappaz2011}.  The numerical velocity from solving \eqref{eq:glenstokes:weak} over $\Lambda(b_h,s_h)$ is denoted $\bu_h$, and its surface trace is denoted $\bu_h|_{s_h}$.  (Note that $\bu_h$ is the numerical, and \emph{not} exact, solution of that Stokes problem.)  We make no assumptions on the solution process for \eqref{eq:glenstokes:weak}, except that it must generate $\bu_h|_{s_h}$.

These preliminaries allow us to define the numerical operator approximating $F^\eps_{\Delta t}$, and the corresponding VI problem.  Following \eqref{eq:definePhi}, \eqref{eq:phisia}, \eqref{eq:regularizedPhi}, and \eqref{eq:regularizedF}, with parameters $\Delta t>0$, $0<\eps\ll 1$, $\eta_0>0$, and for $s_h \in \cK_h$ and $\omega\in\cX$ we define
\begin{align}
F^{\eps,h}_{\Delta t}(s_h)[\omega] &= -\Delta t (1-\eps)\,\int_\Omega \left(\bu_h|_{s_h} \cdot \bn_{s_h}\right) \omega  \label{eq:fe:regularizedF} \\
&\qquad + \Delta t \,\eps\,\Gamma \eta_0^5 \, \int_\Omega |\grad s_h|^2 \grad s_h \cdot \grad \omega + \int_\Omega s_h \omega. \notag
\end{align}
This numerical formula uses the FE function $s_h$ in two ways.  First directly, in the derivatives $\grad s_h$ and $\bn_{s_h}$, and then indirectly to define the Stokes problem domain $\Lambda(b_h,s_h)$, from which $\bu_h$ is found numerically.  Corresponding to VI \eqref{eq:regularizedvi}, with source term $\ell^n$ defined as in \eqref{eq:be:source}, we seek $s_h\in\cK_h$ solving
\begin{equation}
F^{\eps,h}_{\Delta t}(s_h)[\sigma_h-s_h] \ge \ell^n[\sigma_h-s_h] \quad \text{for all } \sigma_h \in \cK_h. \label{eq:fe:regularizedvi}
\end{equation}

Collecting everything above, the following assumptions regarding the data and function spaces will be used when solving the numerical problem \eqref{eq:fe:regularizedvi}, to approximate the exact solution of the regularized implicit time-step problem \eqref{eq:regularizedvi}.

\begin{stdass}
The following data are given:
\renewcommand{\labelenumi}{\arabic{enumi}.}
\begin{enumerate}
\item A bounded, polygonal domain $\Omega\subset\RR^2$.
\item A time-dependent SMB function $a\in C([0,T]; L^{4/3}(\Omega))$.
\item A bed elevation function $b \in C^1(\bar\Omega)$.
\end{enumerate}
We make these continuum definitions:
\begin{enumerate}
\setcounter{enumi}{3}
\item $\cX = W^{1,4}(\Omega)$, with the norm defined in \eqref{eq:norm:Omega}.
\item $\cK = \{\sigma\in\cX\,:\,\sigma|_{\partial \Omega} = b|_{\partial \Omega} \text{ and } \sigma \ge b\}$.
\end{enumerate}
We assume this FE configuration:
\begin{enumerate}
\setcounter{enumi}{5}
\item $\cT_h$ is a triangulation of $\bar\Omega$.
\item $\cX_h \subset \cX$ is the continuous $P_1$ space over $\cT_h$. \label{item:pone}
\item $b_h\in\cX_h$ satisfies $b_h\ge b$ in $\Omega$ and $b_h=b$ along $\partial \Omega$. \label{item:goodbh}
\item $\cK_h = \{\sigma_h\in\cX_h\,:\,\sigma_h|_{\partial \Omega} = b_h|_{\partial \Omega} \text{ and } \sigma_h \ge b_h\}$. \label{item:defineKh}
\end{enumerate}
\end{stdass}

Now we adapt Theorem \ref{thm:abstractestimate}.  Because $s_h\in\cK_h\subset \cK$ by Assumption \ref{item:goodbh}, the first term in \eqref{eq:abstractestimate} is nonpositive, so we replace it by zero, and the technical superset assumption is bypassed.  For $\qq=4$ we have the following result.

\begin{lemma} \label{lem:preglacierapp}  Make the Standard Assumptions, and assume Conjecture \ref{conj:lipschitz} holds for $\rr=4$, and choose parameters $\eps\in(0,1)$, $\eta_0>0$, and $\alpha>0$ so that Conjecture \ref{conj:regcoercive} holds.  Suppose that $s^{n-1}\in\cK$ and define $\ell^n \in \cX'$ by \eqref{eq:be:source}.  For $\Delta t>0$ let $s\in\cK$ be the unique surface elevation satisfying the implicit time-step VI problem \eqref{eq:regularizedvi}, and assume that $s_h\in\cK_h$ solves the numerical problem \eqref{eq:fe:regularizedvi}.  Let $R_h=\max\{\|s\|_\cX,\|s_h\|_\cX\}$, and let $L_\eps(R)$ denote the Lipschitz constant of $\Phi^\eps$ on $B_R\cap\cK$.  Then, with $e_h = s_h-s$ and
\begin{equation}
C_\alpha = \tfrac{3}{4} \left(2\alpha\right)^{-1/3} L_\eps(R_h)^{4/3}, \label{eq:glacierconstant}
\end{equation}
we have
\begin{align}
\frac{\alpha \Delta t}{2} \|e_h\|_\cX^4 + \frac{1}{2} \|e_h\|_{L^2}^2 &\le \inf_{\sigma_h\in\cK_h} \bigg\{\left(F^\eps_{\Delta t}(s)-\ell^n\right)[\sigma_h-s] + \delta_h[s_h-\sigma_h] \label{eq:preglacierestimate} \\
&\hspace{17mm} + \Delta t\, C_\alpha \|\sigma_h - s\|_{\cX}^{4/3} + \frac{1}{2} \|\sigma_h - s\|_{L^2}^2\bigg\}, \notag
\end{align}
where $\delta_h = F^\eps_{\Delta t}(s_h)-F^{\eps,h}_{\Delta t}(s_h) \in \cX'$ is a consistency error at $s_h$.
\end{lemma}

\begin{proof}  Follow the proof of Theorem \ref{thm:abstractestimate}, with $f=F^\eps_{\Delta t}$ and $\qq=4$, but use the coercivity bound $(F^\eps_{\Delta t}(\sigma)-F^\eps_{\Delta t}(s))[\sigma-s] \ge \alpha\Delta t\|\sigma-s\|_\cX^4 + \|\sigma-s\|_{L^2}^2$, and split the Lipschitz term \eqref{eq:abstract:three} as $\Delta t\,(\Phi^\eps(s)-\Phi^\eps(s_h))[s-\sigma_h] + \int_\Omega (s-s_h)(s-\sigma_h)$.  The first part is bounded as in Theorem \ref{thm:abstractestimate}, with $\alpha\Delta t$ and $\Delta t\,L_\eps(R_h)$ in place of $\alpha$ and $L(R_h)$, giving the constant $\Delta t\,C_\alpha$.  The second part is at most $\frac12\|s-s_h\|_{L^2}^2 + \frac12\|s-\sigma_h\|_{L^2}^2$ by Cauchy-Schwarz and Young's inequality.  Because $\cK_h\subset\cK$, the first term in \eqref{eq:abstractestimate} is nonpositive and is dropped.
\end{proof}

Note that the consistency error $\delta_h$ compares the outputs of a continuum operator and a numerical operator at $s_h \in \cK_h$.  However, evaluating $F^\eps_{\Delta t}(s_h)$ requires solving the Stokes problem over the domain $\Lambda(b,s_h)$, with the exact bed, while $F^{\eps,h}_{\Delta t}(s_h)$ involves solving over $\Lambda(b_h,s_h)$, as stated above.  It follows that part of the consistency error is in bed approximation, while part is in other Stokes solver inaccuracies.

The next theorem replaces $\sigma_h$ with an interpolant of the exact solution $s$, and then the meaning of each term in bound \eqref{eq:preglacierestimate} becomes clear.  For $\omega\in\cX = W^{1,4}(\Omega)$ define $\Pi_h(\omega)\in\cX_h$ by
\begin{equation}
\Pi_h(\omega)(x_j) = \max \,\{b_h(x_j), \omega(x_j)\}, \label{eq:definePi}
\end{equation}
over nodes $x_j \in \cT_h$; thus $\Pi_h$ both interpolates and truncates.  If $\sigma\in \cK$ then $\Pi_h(\sigma) \in \cK_h$ because nodal admissibility implies admissibility under Assumptions \ref{item:pone} and \ref{item:goodbh}.

Because $s$ solves \eqref{eq:regularizedvi}, as in Section \ref{sec:abstractestimate} the residual $F^\eps_{\Delta t}(s)-\ell^n \in \cX'$ is represented by a positive Borel measure, supported in the exact (closed) active set $\Omega_A(s)= \left\{x\in\Omega\,:\,s(x)=b(x)\right\}$, which is the ice-free region at $t_n$.  Since $F^\eps_{\Delta t}(s)-\ell^n = \Delta t\,\big((s-s^{n-1})/\Delta t + \Phi^\eps(s) - a^n\big)$, we write this measure as $\Delta t\,\nu_s$, where the positive measure $\nu_s$ represents the regularized, backward-Euler approximation of the residual in \eqref{eq:ncp:residualpos}.

\begin{theorem} \label{thm:glacierapp}  Make the assumptions in Lemma \ref{lem:preglacierapp}, including that Conjectures \ref{conj:lipschitz} and \ref{conj:regcoercive} hold.  Define this scaled surface-trace error from the numerical Stokes approximation,
\begin{equation}
Z_h = \left(1 + C_\infty [L]^{-1} \|s_h\|_{\cX}\right)\, \Big\|\bu_h|_{s_h} - \bu|_{s_h}\Big\|_{L^{4/3}(\Omega)},
\end{equation}
where $C_\infty$ is a Sobolev constant for $\Omega$ (see proof below).  Then, with $e_h=s_h-s$ and $C_\alpha=O(\alpha^{-1/3})$ given by \eqref{eq:glacierconstant},
\begin{align}
\frac{\alpha\Delta t}{4} \|e_h\|_\cX^4 + \frac{1}{2} \|e_h\|_{L^2}^2 &\le \, \Delta t \int_{\Omega_A(s)} (\Pi_h(s) - b)\,d\nu_s &&\text{\textnormal{[term 1]}} \label{eq:glacierestimate} \\
   &\quad\, + \Delta t\left(\tfrac{3}{4}\alpha^{-1/3} Z_h^{4/3} + Z_h \|\Pi_h(s) - s\|_{\cX}\right) &&\text{\textnormal{[term 2]}} \notag \\
   &\quad\, + \Delta t\, C_\alpha \|\Pi_h(s) - s\|_\cX^{4/3} + \frac{1}{2} \|\Pi_h(s) - s\|_{L^2}^2. &&\text{\textnormal{[term 3]}} \notag
\end{align}
\end{theorem}

\begin{proof}  For the first term in \eqref{eq:preglacierestimate} we may write $\Delta t \int_{\Omega_A(s)} (\sigma_h - s) \,d\nu_s$, but observe that $s=b$ in the exact active set.  The consistency error $\delta_h$ in the second term resolves to a difference of velocity traces from \eqref{eq:regularizedF} and \eqref{eq:fe:regularizedF}.  That is, for any $r_h\in \cX_h$,
\begin{equation}
\delta_h[r_h] = \left(F^\eps_{\Delta t}(s_h)-F^{\eps,h}_{\Delta t}(s_h)\right)[r_h] = \Delta t (1-\eps) \int_\Omega \left(\bu_h|_{s_h} - \bu|_{s_h}\right) \cdot \bn_{s_h} \,r_h.
\end{equation}
We estimate the dual norm using H\"older's and Sobolev's inequalities:
\begin{align}
\big|\delta_h[r_h]\big| &\le \Delta t \int_\Omega \Big|\bu_h|_{s_h} - \bu|_{s_h}\Big| \left(1 + |\grad s_h|\right) |r_h| \\
  &\le \Delta t \big\|\bu_h|_{s_h} - \bu|_{s_h}\big\|_{L^{4/3}(\Omega)} \left(1 + C_\infty [L]^{-1} \|s_h\|_{\cX}\right) \|r_h\|_{\cX}. \notag
\end{align}
Defining $\|\delta_h\|_{\ast} = \sup_{0\ne r_h\in\cX_h} \big|\delta_h[r_h]\big| / \|r_h\|_{\cX}$, we have
\begin{equation}
\|\delta_h\|_{\ast} \le \Delta t\, (1 + C_\infty [L]^{-1} \|s_h\|_{\cX}) \big\|\bu_h|_{s_h} - \bu|_{s_h}\big\|_{L^{4/3}(\Omega)} = \Delta t\, Z_h,
\end{equation}
where $C_\infty$ is the Sobolev constant for $\cX\hookrightarrow L^\infty(\Omega)$.  Thus for the second term in \eqref{eq:preglacierestimate} we can use the triangle inequality to write
\begin{equation}
|\delta_h[s_h-\sigma_h]| \le \|\delta_h\|_{\ast} \|s_h-\sigma_h\|_{\cX} \le \|\delta_h\|_{\ast} \|s_h-s\|_{\cX} + \|\delta_h\|_{\ast} \|s-\sigma_h\|_{\cX},
\end{equation}
and thus
\begin{equation}
|\delta_h[s_h-\sigma_h]| \le \frac{\alpha\Delta t}{4} \|s_h-s\|_{\cX}^4 + \tfrac{3}{4}(\alpha\Delta t)^{-1/3} \|\delta_h\|_{\ast}^{4/3} + \|\delta_h\|_{\ast} \|s-\sigma_h\|_{\cX} \label{eq:fromyoungs}
\end{equation}
by Young's inequality with exponents $4$ and $4/3$.  Move the first part to the left side.  Noting that \eqref{eq:preglacierestimate} holds for each fixed $\sigma_h\in\cK_h$, set $\sigma_h=\Pi_h(s)$ to get \eqref{eq:glacierestimate}.
\end{proof}

\section{Discussion of the surface elevation error bound} \label{sec:discussion}

This section addresses the meaning of the Theorem \ref{thm:glacierapp} bound, first by looking at terms separately, then at how the terms scale with long time steps ($\Delta t\to \infty$), and finally regarding the prospect of a convergence proof for coupled FE surface-Stokes solvers ($\Delta t\to 0$).

\subsection{Separate meaning of each term}  For the terms on the right side of \eqref{eq:glacierestimate} we have these comments:

\medskip
\begin{itemize}
\item[term 1:]  This error term comes from margin shape and bed elevation approximation errors in the exact ice-free area, the active set $\Omega_A(s)$ for the unknown surface elevation $s$.  In the unregularized case ($\eps=0$ in VI \eqref{eq:regularizedvi}) the measure $\nu_s$ has density $(b-s^{n-1})/\Delta t - a^n \le -a^n$, so it is bounded by the ablation magnitude over $\Omega_A(s)$.  (See the argument leading to \eqref{eq:be:viearly}; note $\bu|_{s}=\bzero$ and $s=b$ on $\Omega_A(s)$.)

\quad An element with nodes exclusively in the exact ice-free area $\Omega_A(s)$ will have $\Pi_h(s)=b_h$, so the term measures the bed error $b_h-b$ there.  Practically speaking, in those (ablating) portions of the \emph{computed} ice-free area $\Omega_A(s_h)$ which are far from the nearest flowing glacier, by accepting $b_h$ as representing bare land $b$ a modeler may declare that no ice geometry error is made.  However, one cannot neglect term 1 near the unknown glacier margin $\Omega \cap \partial \Omega_A(s)$, the free boundary.  If the numerical bed topography $b_h\approx b$ is inaccurate in a neighborhood of such a location then an ice surface elevation error is inherent, however good the numerical solver.

\medskip
\item[term 2:]  This consistency term quantifies how velocity error from numerically solving the Stokes problem \eqref{eq:glenstokes:weak} over the \emph{numerical} domain $\Lambda(b_h,s_h)$ affects the numerical surface elevation error.  The velocity error is relative to a hypothetical exact solve over $\Lambda(b,s_h)$, i.e.~using the exact bed $b$.  Note that only the velocity traces along $s_h$ are compared.  If $b_h=b$ were to hold and if \eqref{eq:glenstokes:weak} were solved exactly, then this term would be zero.

\quad The factor $Z_h$ scales with the solution norm $\|s_h\|_{\cX}$, and thus a stability bound is needed to understand the right side of \eqref{eq:glacierestimate} entirely in terms of consistency and interpolation errors.  Such a bound is not addressed here.

\quad Note that the two expressions in term 2 scale with different powers of $Z_h$, but the second expression has an interpolation error factor, so it is apparently higher-order under mesh refinement.

\medskip
\item[term 3:]  An interpolation error term like this arises in any C\'ea-type quasi-optimality bound for FE methods on PDEs \cite{Ciarlet2002}.  Separate dependence on the $\cX$ and $L^2$ norms, here and on the left side of \eqref{eq:glacierestimate}, arises because the backward-Euler operator $F^\eps_{\Delta t}$ includes the mass term $\int_\Omega s\omega$; see Definition \eqref{eq:regularizedF}.  Note that $\Pi_h(s)$ is not simply the $P_1$ interpolant of $s$; it includes (nonlinear) truncation into $\cK_h$ to attain discrete admissibility; see Definition \eqref{eq:definePi}.
\end{itemize}

\medskip
\subsection{Steady-state limit}  Suppose that the SMB function $a(x)$ is independent of time.  Dividing VI problem \eqref{eq:regularizedvi} by $\Delta t$ and taking the formal limit $\Delta t\to \infty$ yields the steady-state problem $\Phi^\eps(s)[\sigma-s] \ge \int_\Omega a\,(\sigma-s)$ for all $\sigma\in\cK$.  Conjectures \ref{conj:lipschitz} and \ref{conj:regcoercive} make that problem well-posed by the same argument as in Theorem \ref{thm:regularizedwellposed}.

In this formal steady-state limit \eqref{eq:glacierestimate} says
\begin{align}
  \frac{\alpha}{4} \|e_h\|_\cX^4 &\le \int_{\Omega_A(s)}
  (\Pi_h(s) - b)\,d\nu_s + \tfrac{3}{4}\alpha^{-1/3} Z_h^{4/3} + Z_h \|\Pi_h(s) - s\|_{\cX}
  \label{eq:steadyestimate} \\
     &\quad\, + C_\alpha \|\Pi_h(s) - s\|_\cX^{4/3}, \notag
\end{align}
losing the two $L^2$ mass-term contributions.  (The constant $C_\alpha$ and the factor $Z_h$ do not depend on $\Delta t$.)  Because the source satisfies $\ell^n/\Delta t \to a^n=a$, the previous/initial surface $s^{n-1}$ is forgotten.  The measure $\nu_s$ becomes the steady residual $\Phi^\eps(s)-a \ge 0$, again supported in the ice-free set $\Omega_A(s)$.  In the $\eps\to 0$ unregularized limit the density against Lebesgue measure exists, and it is $-a$ alone in $\Omega_A(s)$.

Thus, the formal steady limit of \eqref{eq:glacierestimate} is the quasi-optimal bound \eqref{eq:steadyestimate} on numerical glacier shape, resembling the steady SIA version in \cite[Theorem 4.3]{JouvetBueler2012}.  Error control in \eqref{eq:steadyestimate} comes only from the $4$-coercivity constant $\alpha$, which depends on the regularization parameters $\eps,\eta_0$.  In contrast, when $\Delta t$ is finite, the $L^2$ mass term in \eqref{eq:glacierestimate} gives additional control, making the time-stepping problem better conditioned than the steady one.

\subsection{Numerical convergence}  Theorem \ref{thm:glacierapp} is not a convergence theorem.  It bounds the error of a single time step, it is conditional on Conjectures \ref{conj:lipschitz} and \ref{conj:regcoercive}, and it says nothing about how errors accumulate over many steps.  The following remarks therefore involve some speculation.

Regarding temporal refinement, as $\Delta t\to 0$ over a fixed mesh, terms 1 and 2 in \eqref{eq:glacierestimate}, and the $\cX$-norm part of term 3, are $O(\Delta t)$, and the bound degenerates to $\|e_h\|_{L^2}\lesssim\|\Pi_h(s)-s\|_{L^2}$.  Simply summing such $O(1)$ single-step errors would not converge.  Instead, space-time convergence needs stability of the time-step map with respect to $s^{n-1}$, though pursuing this is outside of our scope.

Consider 2D mesh refinement $h\to 0$ for a fixed time step.  Each term in \eqref{eq:glacierestimate} should shrink: term 1 through more accurate bed interpolation $b_h\approx b$ near the free boundary, term 2 through a more accurate Stokes solution (presumably requiring 3D refinement), and term 3 through better approximation of $s$ in $\cX_h$.

Precise convergence rates are harder.  If $s\in W^{2,\mu}(\Omega)$ for some $\mu\in[4,\infty]$, and ignoring the truncation occurring in $\Pi_h$, interpolation theory \cite[Theorem 3.1.6]{Ciarlet2002} gives $\|\Pi_h(s)-s\|_\cX = O(h)$, so term 3 is $O(h^{4/3})$ and one might expect $\|e_h\|_\cX = O(h^{1/3})$.  (Theorem 4.3 in \cite{JouvetBueler2012} makes a comparable assumption for an SIA model.)  However, $s\in W^{2,\mu}$ would force the surface to be tangent to the bed at the glacier margin, which real glaciers do not do; $s\in W^{1,4}$ is credible but $s\in W^{2,\mu}$ is not.

Considering the Stokes domain $\Lambda=\Lambda(b_h,s_h)$ used in term 2 as fixed for the moment, and assuming regularity for $\bu,p$, and a discrete inf-sup condition for the elements, mixed FE theory gives velocity errors $O(h^{\kappa/2})$ for some $\pp \le \kappa \le 2$ \cite{JouvetRappaz2011}.  However, the constant depends on the domain $\Lambda$, equivalently on $s_h$, which itself changes under refinement.  A convergence proof would need uniform control on all $h$-dependent constants, and we do not currently know how to obtain it.

\section{Conclusion} \label{sec:conclusion}

The major result of this paper is Theorem \ref{thm:glacierapp}.  Though it is subject to nontrivial conjectures about the continuum problem, it provides a bound on the finite element (FE) surface elevation error committed in a single implicit time step of a fully-coupled and non-shallow glacier dynamics model.

In developing this result we regard the velocity trace from a Stokes stress balance, more precisely the surface-motion term in the free-surface kinematical equation \eqref{eq:ske}, as the output of an abstract map acting on the surface elevation.  The implicit time-step problem is then a regularized variational inequality  \eqref{eq:regularizedvi}, the solution of which simultaneously determines the map-plane extent of glaciation, the surface elevation function itself, and the velocity and pressure within the 3D ice domain (determined by that solution surface elevation).  Motivated by the shallow ice approximation, the regularization adds a small amount of diffusion.

We know of no existing literature which is comparably-rigorous to what is done here for Stokes-based glacier evolution.  On the one hand, numerical solvers for this kind of problem are plentiful \cite{AhlkronaLofgrenHenry2025, Brinkerhoff2023, Chengetal2020, IsaacStadlerGhattas2015, Jouvetetal2008, Lengetal2012, WirbelJarosch2020}.  On the other, a candidate function space for the surface elevation (or thickness) must obviously be proposed in any theory of well-posedness, or for any precise FE error bound on the elevation.  Noting that the theory of shallow models \cite{JouvetBueler2012,PiersantiTemam2023} provides some guidance, the current paper seems to be the first one to propose a Sobolev space for the surface elevation in a non-shallow theory, or to include the behavior of the Stokes surface motion in any error analysis.

Related recent work provides an $L^2$ bound on time-dependent surface elevations \cite[Theorem 3.3]{Tominecetal2026} in the coupled surface-Stokes model.  However, as noted by those authors, the $L^2$ norm in this bound must play an incomplete role in the well-posedness theory of the coupled problem.  The evolving Stokes domain must have sufficient boundary regularity so as to generate a well-defined surface motion at each solution, but an $L^2$ bound on the surface elevation function cannot provide that.

We do not possess a proof, or dis-proof, of the well-posedness of regularized problem \eqref{eq:regularizedvi}.  The heart of the matter seems to be that the unregularized surface motion operator from the Stokes model acts $4$-coercively in the $W^{1,4}$ norm on most admissible pairs of surfaces.  Numerical evidence to this effect is in Appendix \ref{app:numerical}.  Further evidence in Appendix \ref{app:numerical} says that a small amount of diffusion, motivated by the shallow ice approximation, appears sufficient to generate a $4$-coercive operator, and this motivates our Conjecture \ref{conj:regcoercive}.

%\clearpage\newpage
\bibliographystyle{siamplain}
\bibliography{estimate}

@book{Acheson1990,
  title={Elementary {F}luid {D}ynamics},
  author={Acheson, DJ},
  year={1990},
  publisher={Oxford University Press}
}

@unpublished{AhlkronaLofgrenHenry2025,
  author={J. Ahlkrona and A. Löfgren and C. Henry},
  title={A stable, fully implicit, second order method for viscous free surface Stokes flow},
  note={submitted},
  year={2025},
}

@book{AscherPetzold1998,
    AUTHOR = {Ascher, U. and L. Petzold},
     TITLE = {Computer {M}ethods for {O}rdinary {D}ifferential {E}quations and
              {D}ifferential-algebraic {E}quations},
 PUBLISHER = {SIAM Press},
   ADDRESS = {Philadelphia},
      YEAR = {1998},
}

@article{Aschwandenetal2012,
  title={An enthalpy formulation for glaciers and ice sheets},
  author={A. Aschwanden and E. Bueler and C. Khroulev and H. Blatter},
  journal={Journal of Glaciology},
  volume={58},
  number={209},
  pages={441--457},
  doi={10.3189/2012JoG11J088},
  year={2012},
}

@article{Belenkietal2012,
author = {Belenki, L. and Berselli, L. C. and Diening, L. and R\r{u}\v{z}i\v{c}ka, M.},
title = {On the finite element approximation of $p$-{S}tokes systems},
journal = {SIAM J. Numer. Anal.},
volume = {50},
number = {2},
pages = {373-397},
year = {2012},
doi = {10.1137/10080436X},
}

@Article{BensonMunson2006,
  AUTHOR    = {S. Benson and T. Munson},
  TITLE     = {Flexible complementarity solvers for large-scale applications},
  JOURNAL   = {Optimization Methods and Software},
  VOLUME    = {21},
  NUMBER    = {1},
  PAGES     = {155--168},
  doi={10.1080/10556780500065382},
  YEAR      = {2006},
}

@Article{Bodvardsson1955,
  AUTHOR    = {G. Bodvardsson},
  TITLE     = {On the flow of ice-sheets and glaciers},
  JOURNAL   = {J{\"o}kull},
  VOLUME    = {5},
  PAGES     = {1--8},
  YEAR      = {1955},
}

@book{BoffiBrezziFortin2013,
  title={{Mixed Finite Element Methods and Applications}},
  author={D. Boffi and F. Brezzi and M. Fortin},
  volume={44},
  year={2013},
  publisher={Springer}
}

@Article{Brinkerhoff2023,
  author={D. J. Brinkerhoff},
  title={Compatible finite elements for glacier modeling},
  journal={Computing in Science \& Engineering},
  year={2023},
  volume={25},
  number={3},
  pages={18--28},
  doi={10.1109/MCSE.2023.3305864},
}

@Article{Bueler2021conservation,
  author   = {E. Bueler},
  title    = {Conservation laws for free-boundary fluid layers},
  journal  = {SIAM J. Appl. Math.},
  volume   = {81},
  number   = {5},
  pages    = {2007--2032},
  year     = {2021},
  doi      = {10.1137/20M135217X},
}

@article{Bueler2023,
  title={Performance analysis of high-resolution ice-sheet simulations},
  author={E. Bueler},
  journal={J. Glaciol.},
  volume={69},
  number={276},
  pages={930--935},
  doi={10.1017/jog.2022.113},
  year={2023}
}

@Article{BuelerFarrell2024,
  author   = {E. Bueler and P. Farrell},
  title    = {A full approximation scheme multilevel method for nonlinear variational inequalities},
  journal={SIAM J. Sci. Comput.},
  volume={46},
  number={4},
  pages={A2421--A2444},
  doi = {10.1137/23M1594200},
  year={2024},
}

@Article{Bueleretal2005,
  author	= {E. Bueler and C. S. Lingle and J. A. Kallen-Brown and D.
		       N. Covey and L. N. Bowman},
  title		= {Exact solutions and verification of numerical models for isothermal
		       ice sheets},
  journal	= {J. Glaciol.},
  volume	= {51},
  number	= {173},
  pages		= {291--306},
  year		= {2005},
  doi       = {10.3189/172756505781829449},
}

@Article{Calvoetal2003,
  author  = {N. Calvo and others},
  title   = {On a doubly nonlinear parabolic obstacle problem modelling
             ice sheet dynamics},
  year    = {2003},
  journal = {SIAM J. Appl. Math.},
  volume  = {63},
  number  = {2},
  pages   = {683--707},
  doi     = {10.1137/S0036139901385345},
}

@Article{Chengetal2020,
AUTHOR = {Cheng, G. and L\"otstedt, P. and von Sydow, L.},
TITLE = {A full {S}tokes subgrid scheme in two dimensions for simulation of grounding line migration in ice sheets using {Elmer/ICE} (v8.3)},
JOURNAL = {Geoscientific Model Development},
VOLUME = {13},
YEAR = {2020},
NUMBER = {5},
PAGES = {2245--2258},
DOI = {10.5194/gmd-13-2245-2020}
}

@Article{ChoeLewis1991,
  author={H. J. Choe and J. L. Lewis},
  title={On the obstacle problem for quasilinear elliptic equations of $p$-{L}aplacian type},
  journal={SIAM J. Math. Anal.},
  volume={22},
  number={3},
  pages={623--638},
  year={1991},
}

@Article{Chow1989,
AUTHOR = {Chow, S. S.},
TITLE = {Finite element error estimates for non-linear elliptic equations of monotone type},
JOURNAL = {Numer. Math.},
VOLUME = {54},
YEAR = {1989},
PAGES = {373--393},
DOI = {10.1007/BF01396320}
}

@Book{Ciarlet2002,
  AUTHOR    = {P. Ciarlet},
  TITLE     = {The {F}inite {E}lement {M}ethod for {E}lliptic {P}roblems},
  NOTE      = {Reprint of the 1978 original},
  PUBLISHER = {SIAM Press},
  ADDRESS   = {Philadelphia},
  YEAR      = {2002},
}

@TechReport{Cogleyetal2011,
  Author = {J.G. Cogley and R. Hock and L.A. Rasmussen and A.A. Arendt and A. Bauder and R.J. Braithwaite and P. Jansson and G. Kaser and M. Möller and L. Nicholson and M. Zemp},
  Title = {Glossary of Glacier Mass Balance and Related Terms},
  type = {IHP-VII Technical Documents in Hydrology},
  number = {86},
  institution = {UNESCO-IHP},
  address = {Paris},
  Year = {2011},
}

@article{Daners2003,
  author={Daners, D.},
  title={Dirichlet problems on varying domains},
  journal={Journal of Differential Equations},
  volume={188},
  number={2},
  pages={591--624},
  year={2003},
}

@article{Durandetal2009,
author = {Durand, G. and Gagliardini, O. and de Fleurian, B. and Zwinger, T. and Le Meur, E.},
title = {Marine ice sheet dynamics: Hysteresis and neutral equilibrium},
journal = {Journal of Geophysical Research: Earth Surface},
volume = {114},
number = {F3},
doi = {https://doi.org/10.1029/2008JF001170},
year = {2009},
}

@Article{EchelmeyerKamb1986,
  title={Stress-gradient coupling in glacier flow: {II}. {Longitudinal} averaging in the flow response to small perturbations in ice thickness and surface slope},
  author={K. A. Echelmeyer and B. Kamb},
  journal={J. Glaciol.},
  volume={32},
  number={111},
  pages={285--298},
  doi={10.3189/S0022143000015616},
  year={1986},
}

@Book{Elmanetal2014,
  AUTHOR    = {H. C. Elman and D. J. Silvester and A. J. Wathen},
  TITLE     = {Finite {E}lements and {F}ast {I}terative {S}olvers: with
               {A}pplications in {I}ncompressible {F}luid {D}ynamics},
  EDITION   = {2nd},
  PUBLISHER = {Oxford University Press},
  ADDRESS   = {Oxford, UK},
  YEAR      = {2014},
}

@Book{Evans2010,
  AUTHOR    = {L. C. Evans},
  TITLE     = {Partial {D}ifferential {E}quations},
  SERIES    = {Graduate Studies in Mathematics},
  PUBLISHER = {American Mathematical Society},
  ADDRESS   = {Providence},
  EDITION   = {2nd},
  YEAR      = {2010},
}

@Book{FacchineiPang2003,
  author={Facchinei, F. and Pang, J.-S.},
  title={{Finite-Dimensional Variational Inequalities and Complementarity Problems}},
  volume={1},
  year={2003},
  publisher={Springer}
}

@article{Falk1974,
  author={R. S. Falk},
  title={Error estimates for the approximation of a class of variational inequalities},
  journal={Mathematics of Computation},
  volume={28},
  number={128},
  pages={963--971},
  year={1974}
}

@Book{Grisvard1985,
  AUTHOR    = {P. Grisvard},
  TITLE     = {{Elliptic Problems in Nonsmooth Domains}},
  PUBLISHER = {Pitman},
  ADDRESS   = {Boston},
  YEAR      = {1985},
}

@Article{GoldsbyKohlstedt2001,
    AUTHOR = {D. L. Goldsby and D. L. Kohlstedt},
     TITLE = {Superplastic deformation of ice: experimental observations},
   JOURNAL = {J. Geophys. Res.},
    VOLUME = {106},
    NUMBER = {B6},
     PAGES = {11017--11030},
     doi={10.1029/2000JB900336},
      YEAR = {2001},
}

@Book{GreveBlatter2009,
  AUTHOR    = {R. Greve and H. Blatter},
  TITLE     = {Dynamics of {I}ce {S}heets and {G}laciers},
  SERIES    = {Advances in Geophysical and Environmental Mechanics and
               Mathematics},
  PUBLISHER = {Springer},
  ADDRESS   = {Berlin, Germany},
  YEAR      = {2009},
}

@Article{Halfar1981,
  author={P. Halfar},
  title={On the dynamics of the ice sheets},
  journal={J. Geophys. Res.},
  volume={86},
  number={C11},
  pages={11065--11072},
  doi={10.1029/JC086iC11p11065},
  year={1981},
}

@manual{Hametal2023,
  title        = {Firedrake User Manual},
  author       = {David A. Ham and Paul H. J. Kelly and Lawrence Mitchell and Colin J. Cotter and Robert C. Kirby and Koki Sagiyama and Nacime Bouziani and Sophia Vorderwuelbecke and Thomas J. Gregory and Jack Betteridge and Daniel R. Shapero and Reuben W. Nixon-Hill and Connor J. Ward and Patrick E. Farrell and Pablo D. Brubeck and India Marsden and Thomas H. Gibson and Miklós Homolya and Tianjiao Sun and Andrew T. T. McRae and Fabio Luporini and Alastair Gregory and Michael Lange and Simon W. Funke and Florian Rathgeber and Gheorghe-Teodor Bercea and Graham R. Markall},
  organization = {Imperial College London and University of Oxford and Baylor University and University of Washington},
  edition      = {First edition},
  year         = {2023},
  month        = {5},
  doi          = {10.25561/104839},
}

@Article{IsaacStadlerGhattas2015,
  title={Solution of nonlinear {S}tokes equations discretized by high-order finite elements on nonconforming and anisotropic meshes, with application to ice sheet dynamics},
  author={Isaac, T. and Stadler, G. and Ghattas, O.},
  journal={SIAM J. Sci. Comput.},
  volume={37},
  number={6},
  pages={B804--B833},
  doi={10.1137/140974407},
  year={2015},
}

@Article{Jouvetetal2008,
  title={A new algorithm to simulate the dynamics of a glacier: theory and applications},
  author={Jouvet, G. and Picasso, M. and Rappaz, J. and Blatter, H.},
  journal={J. Glaciol.},
  volume={54},
  number={188},
  pages={801--811},
  doi={10.3189/002214308787780049},
  year={2008},
}

@Article{Jouvetetal2011,
  title={Existence and stability of steady-state solutions of the shallow-ice-sheet equation by an energy-minimization approach},
  author={Jouvet, G. and Rappaz, J. and Bueler, E. and Blatter, H.},
  journal={J. Glaciol.},
  volume={57},
  number={202},
  pages={345--354},
  year={2011},
}

@Article{JouvetBueler2012,
  author	= {G. Jouvet and E. Bueler},
  title		= {Steady, shallow ice sheets as obstacle problems: well-posedness and finite element approximation},
  journal	= {SIAM J. Appl. Math.},
  volume	= {72},
  number	= {4},
  pages		= {1292--1314},
  doi={10.1137/110856654},
  year		= {2012},
}

@Article{JouvetRappaz2011,
  title={Analysis and finite element approximation of a nonlinear stationary {S}tokes problem arising in glaciology},
  author={Jouvet, G. and Rappaz, J.},
  journal={Advances in Numerical Analysis},
  doi={10.1155/2011/164581},
  year={2011},
}

@Book{KinderlehrerStampacchia1980,
  AUTHOR    = {D. Kinderlehrer and G. Stampacchia},
  TITLE     = {An {I}ntroduction to {V}ariational {I}nequalities and
               their {A}pplications},
  PUBLISHER = {Academic Press},
  ADDRESS   = {New York},
  YEAR      = {1980},
}

@article{Lengetal2012,
author = {Leng, Wei and Ju, Lili and Gunzburger, Max and Price, Stephen and Ringler, Todd},
title = {A parallel high-order accurate finite element nonlinear {S}tokes ice sheet model and benchmark experiments},
journal = {J. Geophys. Res.: Earth Surface},
volume = {117},
number = {F1},
doi = {https://doi.org/10.1029/2011JF001962},
year = {2012},
}

@Book{LiebLoss1997,
  AUTHOR    = {E. H. Lieb and M. Loss},
  TITLE     = {Analysis},
  PUBLISHER = {American Mathematical Society},
  ADDRESS   = {Providence},
  SERIES    = {Graduate Studies in Mathematics},
  VOLUME    = {14},
  YEAR      = {1997},
}

@article{LofgrenAhlkronaHelanow2022,
  author={A. L{\"o}fgren and J. Ahlkrona and C. Helanow},
  title={Increasing stable time-step sizes of the free-surface problem arising in ice-sheet simulations},
  journal={Journal of Computational Physics: X},
  volume={16},
  number={100114},
  year={2022},
  doi={10.1016/j.jcpx.2022.100114},
}

@article{Minty1963,
author = {G. J. Minty },
title = {On a ``monotonicity'' method for the solution of nonlinear equations in {B}anach spaces},
journal = {Proceedings of the National Academy of Sciences},
volume = {50},
number = {6},
pages = {1038-1041},
year = {1963},
doi = {10.1073/pnas.50.6.1038},
}

@Article{Morlighemetal2017,
author = {Morlighem, M. and others},
title = {{BedMachine v3}: {C}omplete bed topography and ocean bathymetry mapping of {G}reenland from multibeam echo sounding combined with mass conservation},
journal = {Geophysical Research Letters},
volume = {44},
number = {21},
pages = {11051--11061},
doi = {10.1002/2017GL074954},
year = {2017},
}

@article{Pattynetal2008,
  title={Benchmark experiments for higher-order and full-{S}tokes ice sheet models {(ISMIP--HOM)}},
  author={Pattyn, F. and Perichon, L. and others},
  journal={The Cryosphere},
  volume={2},
  number={2},
  pages={95--108},
  doi={10.5194/tc-2-95-2008},
  year={2008},
}

@article{PiersantiTemam2023,
  author={P. Piersanti and R. Temam},
  title={On the dynamics of grounded shallow ice sheets: Modeling and analysis},
  journal={Advances in Nonlinear Analysis},
  volume={12},
  year={2023},
  doi={10.1515/anona-2022-0280},
}

@article{Pompe2003,
  title={Korn's first inequality with variable coefficients and its generalization},
  author={W. Pompe},
  journal={Commentationes Mathematicae Universitatis Carolinae},
  volume={44},
  number={1},
  pages={57--70},
  year={2003},
}

@Article{PralongFunk2005,
author = {Pralong, A. and Funk, M.},
title = {Dynamic damage model of crevasse opening and application to glacier calving},
journal = {J. Geophys. Res.: Solid Earth},
volume = {110},
number = {B1},
doi = {https://doi.org/10.1029/2004JB003104},
year = {2005},
}

@Article {SchoofHewitt2013,
  AUTHOR    = {C. Schoof and I. J. Hewitt},
  TITLE     = {Ice-sheet dynamics},
  JOURNAL   = {Annual Review of Fluid Mechanics},
  VOLUME    = {45},
  PAGES     = {217--239},
  YEAR      = {2013},
  doi       = {10.1002/jgrf.20146},
}

@article{Tominecetal2026,
  title={Free-surface {S}tokes problem: stability estimates and time-step improvements},
  author={Tominec, Igor and Lundgren, Lukas and L{\"o}fgren, Andr{\'e} and Ahlkrona, Josefin},
  journal={IMA Journal of Numerical Analysis},
  year={2026},
  doi={10.1093/imanum/drag117},
}

@Article{Winkelmannetal2011,
AUTHOR = {Winkelmann, R. and Martin, M. A. and Haseloff, M. and Albrecht, T. and Bueler, E. and Khroulev, C. and Levermann, A.},
TITLE = {The {P}otsdam {P}arallel {I}ce {S}heet {M}odel ({PISM-PIK}) {P}art 1: {M}odel description},
JOURNAL = {The Cryosphere},
VOLUME = {5},
PAGES = {715--726},
DOI={10.5194/tc-5-715-2011},
YEAR = {2011},
}

@Article {WirbelJarosch2020,
  AUTHOR    = {A Wirbel and A Jarosch},
  TITLE     = {Inequality-constrained free-surface evolution in a full {S}tokes ice flow model (evolve\_glacier v1.1)},
  JOURNAL   = {Geoscientific Model Development},
  VOLUME    = {13},
  NUMBER    = {12},
  PAGES     = {6425--6445},
  YEAR      = {2020},
  doi       = {10.5194/gmd-13-6425-2020},
}

\appendix
\section{Proof of Lemma \ref{lem:stokesapriori}} \label{app:provestokesapriori}

The needed Poincar\'e and Korn inequalities are stated next.  Recall that $\Lambda$ is assumed to be a bounded Lipschitz domain (Section \ref{sec:stokes}).  Inequality \eqref{eq:poincare} is proven below using the form \eqref{eq:icydomain} of $\Lambda$ and the zero trace of $\bv\in\cV$ along $\Gamma_b$.  Inequality \eqref{eq:korns} follows by setting $F(x)$ to the identity in Corollary 4.1 of \cite{Pompe2003}, which applies because $\Gamma_b$ has positive measure.

\begin{lemma} \label{lem:poincarekorns}
There exist dimensionless constants $C$, depending on $\pp$, $[H]$, and the geometry of $\Lambda$, so that for all $\bv \in \cV$,
\begin{equation}
\int_\Lambda |\bv|^\pp \le C [H]^\pp \int_\Lambda |\grad\bv|^\pp, \label{eq:poincare}
\end{equation}
thus $\|\bv\|_{\cV}^\pp \le (C + 1) [H]^\pp \int_\Lambda |\grad\bv|^\pp$, and
\begin{equation}
\int_\Lambda |\grad\bv|^\pp \le C \int_\Lambda |D\bv|^\pp. \label{eq:korns}
\end{equation}
\end{lemma}

\begin{proof}[Proof of \eqref{eq:poincare}]
Let $H_{\max} = \sup_x (s(x)-b(x))$ be the maximum thickness of $\Lambda = \{(x,z)\,:\,b(x)<z<s(x)\}$.  By density it suffices to consider smooth $\bv$ which vanish in a neighborhood of $\Gamma_b$.  For $b(x)<z<s(x)$, H\"older's inequality gives
$$|\bv(x,z)| \le \int_{b(x)}^z |\partial_z \bv(x,z')|\,dz' \le (z-b(x))^{1/\pp'} \left(\int_{b(x)}^{s(x)} |\partial_z \bv(x,z')|^\pp\,dz'\right)^{1/\pp}.$$
Raising to the power $\pp$, integrating in $z$ over $(b(x),s(x))$, and then integrating in $x$, yields
$$\int_\Lambda |\bv|^\pp \le \frac{H_{\max}^\pp}{\pp} \int_\Lambda |\partial_z \bv|^\pp \le \frac{H_{\max}^\pp}{\pp} \int_\Lambda |\grad\bv|^\pp.$$
This is \eqref{eq:poincare} with dimensionless constant $C = (H_{\max}/[H])^\pp/\pp$.
\end{proof}

The following proof is hinted-at in \cite{JouvetRappaz2011}; it takes a standard approach.

\begin{proof}[Proof of Lemma \ref{lem:stokesapriori}]
From \eqref{eq:glenstokes:weak} and $\bu \in\cV_0$ it follows that
\begin{equation}
0= F_\Lambda(\bu,p)[\bu,p] = \int_\Lambda 2 \nu(D\bu) D\bu : D\bu - \rhoi \bg \cdot \bu.  \label{eq:stokes:substituteu}
\end{equation}
Apply \eqref{eq:korns}, the facts that $\pp>0$ and $(\pp-2)/2 \le 0$, and equation \eqref{eq:glen}:
\begin{align}
\int_\Lambda |\grad\bu|^\pp &\le C \int_\Lambda |D\bu|^\pp \le C \int_\Lambda \left(|D\bu|^2 + \mu_0\right)^{(\pp-2)/2} \left(D\bu:D\bu + \mu_0\right) \label{eq:stokes:startapriori} \\
	&\le C \left[\mu_0^{\pp/2} |\Lambda| + (2 \nu_\pp)^{-1} \int_\Lambda 2\nu(D\bu) D\bu:D\bu\right]. \notag
\end{align}
By equation \eqref{eq:stokes:substituteu}, Cauchy-Schwarz, H\"older's inequality, and definition \eqref{eq:vnorm}, we have $\int_\Lambda 2\nu(D\bu) D\bu:D\bu = \int_\Lambda \rhoi \bg\cdot\bu \le \rhoi |\bg| |\Lambda|^{1/\pp'} \|\bu\|_\cV$.  Combining this with the consequence $\|\bu\|_\cV^\pp \le (C+1)[H]^\pp \int_\Lambda |\grad\bu|^\pp$ of \eqref{eq:poincare}, we thus have
\begin{equation}
\|\bu\|_\cV^\pp \le C [H]^\pp \left[\mu_0^{\pp/2} |\Lambda| + (2 \nu_\pp)^{-1} \rhoi |\bg| |\Lambda|^{1/\pp'} \|\bu\|_\cV\right]. \label{eq:stokes:workapriori}
\end{equation}
Now let $z=\|\bu\|_\cV$.  We have proved that $z^\pp \le c_0 + c_1 z$ for $\pp>1$ and constants $c_i>0$.  Furthermore $g(y) = y^\pp - c_1 y - c_0$ is continuous and convex on $[0,\infty)$, with $g(0)=-c_0<0$, and $g(y) \to +\infty$ as $y \to +\infty$ since $\pp>1$.  Thus there exists a unique positive root $\tilde y>0$, and $g'(\tilde y) \ge (g(\tilde y)-g(0))/\tilde y > 0$ by convexity, so by the implicit function theorem $\tilde y = f(\pp,c_0,c_1)$ depends continuously on its arguments.  Since $g(z)\le 0$ we have $z \le \tilde y$.  This proves \eqref{eq:stokesapriori} with $C=\tilde y$.
\end{proof}

\section{The surface motion map is not strictly monotone} \label{app:noncoercive}

It is reasonable to suppose that a glacier flow model will produce no flow for an admissible surface elevation $s\in \cX=W^{1,\rr}(\Omega)$ which is horizontal everywhere on the ice.  In terms of the inactive set $\Omega_I=\{x\in\Omega\,:\,s(x)>b(x)\}$ and the domain $\Lambda(b,s)$ from \eqref{eq:domainfroms}, this property says that if $\grad s=0$ a.e.~in $\Omega_I$ then $\bu=\bzero$ on $\Lambda(b,s)$.

This property holds for the non-sliding Stokes model considered in the paper because $\bu=\bzero$ and $p=\rhoi |\bg| (s-z)$ (i.e.~hydrostatic) solve the model, and the solution is unique.  The property also holds for most sliding laws, and in SIA models \cite{JouvetBueler2012}.

\begin{figure}[ht]
\begin{center}
\mbox{\includegraphics[width=0.4\textwidth]{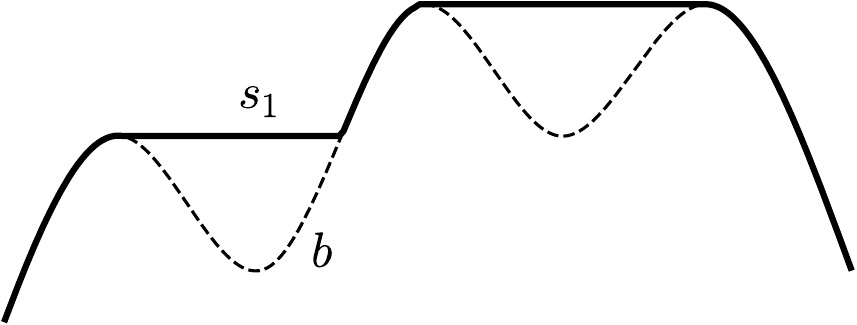} \qquad \includegraphics[width=0.4\textwidth]{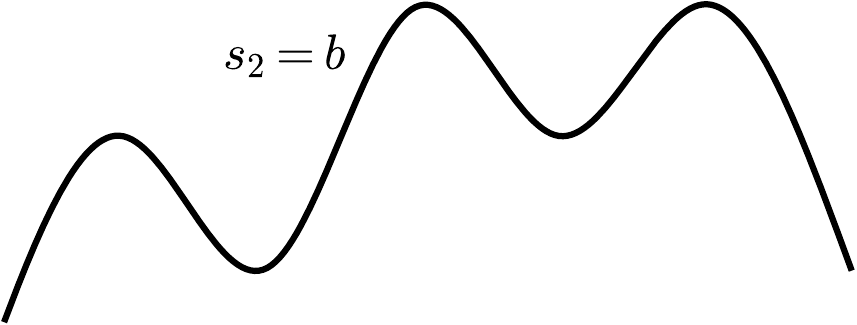}}
\end{center}
\caption{For a smooth bedrock elevation $b$ with strict local minima, one surface elevation state $s_1$ fills in the minima with flat ice (left), while the other is ice-free $s_2=b$ (right).  Because neither state generates flow, $\bu|_{s_i}=\bzero$, the surface motion map $\Phi$ in equation \eqref{eq:definePhi} cannot be strictly monotone.}
\label{fig:noncoercive}
\end{figure}

Using this property, for any bedrock elevation $b$ possessing strict local minima in the interior of $\Omega$ we can construct distinct surface elevations $s_1,s_2\in \cK$, $\|s_1-s_2\|_{\cX}>0$, with zero surface velocity $\bu|_{s_i}=\bzero$ \eqref{eq:defineus}, and therefore zero surface motion from \eqref{eq:definePhi}: $\Phi(s_i)=0$.  The construction (Figure \ref{fig:noncoercive}) simply fills in the local minima with flat ice for surface $s_1$, and leaves $s_2=b$ as bare ground.  Note that $s_1$ can be constructed with any number of local minima filled in, and/or with the minima only partially filled, so generally there are infinitely many surfaces of this type.  Two consequences of this construction are as follows.

\begin{proposition} \label{prop:noncoercive}
The surface motion map $\Phi$ defined in \eqref{eq:definePhi} is not $\qq$-coercive, for any $\qq$, nor is it strictly monotone.  (See Definition \ref{def:monotonecoercive}.)
\end{proposition}

\begin{proof}
Note $(\Phi(s_1) - \Phi(s_2))[s_1-s_2]=0$ while $\|s_1-s_2\|_\cX > 0$.
\end{proof}

\begin{proposition} \label{prop:notunique}
When the SMB is identically zero, $a=0$, there are smooth bedrock elevations for which the steady-state version of the glacier geometry model \eqref{eq:ncp}, or of any other model with the same no-flow property (above), will have more than one solution.
\end{proposition}

\begin{proof}
The strong (NCP) form of the steady-state and $a=0$ model is
	$$s - b \ge 0, \quad - \bu|_s \cdot \bn_s \ge 0, \quad (s - b) \left(- \bu|_s \cdot \bn_s\right) = 0.$$
If $b$ is a smooth function with strict local minima then the above construction of $s_1$ and $s_2$ gives distinct solutions.
\end{proof}

Proposition \ref{prop:notunique} answers negatively, via a surprisingly simple construction which seems not to appear in the literature, the uniqueness question for the general-bed, steady-state SIA, at least for certain elevation-\emph{in}dependent SMB functions.  (The counter-example needs only $a=0$ on the filled basins and $a\le 0$ elsewhere.)  This question has apparently been open since the proof of existence \cite{JouvetBueler2012} for that model, which applies to general $a(x)$.  Non-uniqueness here is very different from the better-known non-uniqueness \cite{Bodvardsson1955} and non-existence \cite{Jouvetetal2011} results for certain elevation-dependent SMB models.

\section{Experiments on coercivity of the surface motion map} \label{app:numerical}

Conjecture \ref{conj:regcoercive} is part of a well-posedness framework for the continuum implicit time-step problem for glaciers, VI problem \eqref{eq:regularizedvi}.  In combination with \emph{a priori} Theorem \ref{thm:abstractestimate}, this leads to a bound on FE surface elevation error (Theorem \ref{thm:glacierapp}).

The reasonableness of the Conjecture can be tested by sampling from numerical simulations.  The experiments here, using Python and the Firedrake FE library \cite{Hametal2023},\footnote{Source code is at the public repository \href{https://github.com/bueler/glacier-fe-estimate}{\texttt{github.com/bueler/glacier-fe-estimate}}, in the \texttt{py/} directory.  The codes call the library at \href{https://github.com/bueler/stokes-extrude}{\texttt{github.com/bueler/stokes-extrude}}.} are not intended to demonstrate any particular time-stepping.  They simply generate admissible surface elevation pairs $\sigma,s\in\cK \subset \cX = W^{1,4}(\Omega)$, to use as samples in a coercivity test.  For each pair of surfaces we evaluate $4$-coercivity ratios:
\begin{equation}
\rho^\eps(\sigma,s) = \frac{\left(\Phi^\eps(\sigma) - \Phi^\eps(s)\right)[\sigma-s]}{\|\sigma-s\|_{\cX}^4}. \label{eq:Phiratio}
\end{equation}
(Compare inequality \eqref{eq:qcoercive}.)  In the experiments we use $\eps\in\{0,0.01\}$ only, and $\eta_0=500$ m in \eqref{eq:fe:regularizedF}.  Appendix \ref{app:noncoercive} shows that for bumpy beds there exist pairs that give $\rho^0(\sigma,s) = 0$.

\begin{figure}[ht]
\centering
\includegraphics[width=0.35\textwidth]{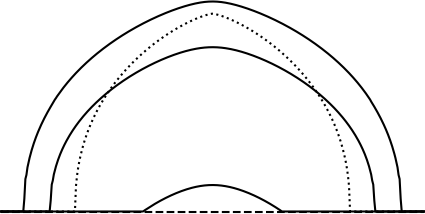} \quad \includegraphics[width=0.37\textwidth]{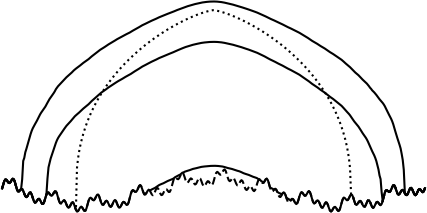}
\caption{Two bed profiles, flat and rough, define constraint sets $\cK_h^i\subset\cX_h$, $i=1,2$, in the numerical experiment.  For each $\cK_h^i$, runs using different constant values of the SMB generated a large number of admissible surface elevations, representing a range of glacier sizes covering different portions of the domain.}
\label{fig:cases}
\end{figure}

If the set of all possible ratios $\{\rho^\eps(\sigma,s)\}_{\sigma,s\in\cK}$ were bounded below by a positive constant $\alpha>0$, then this would verify Conjecture \ref{conj:regcoercive}.  However, an experiment obviously allows only finite sampling using particular data and discretizations.  Here the domain is the 1D interval $\Omega=(-L,L)$ with $L=100$ km, meshed into equal intervals.  The FE space $\cX_h\subset \cX$ is the set of continuous $P_1$ (piecewise-linear) functions, thus the Stokes domain $\Lambda(b_h,s_h)$ is polygonal.  Regarding the data, two bed profiles (Figure \ref{fig:cases}) are considered, \emph{flat} ($b_1(x)=0$) and \emph{rough} ($b_2(x)$); the latter is a linear combination of several sinusoids down to 4 km wavelength.  Three constant SMB values $a_j\in\{-2.5,0.0,1.0\}\times 10^{-7}$\, $\text{m}\,\text{s}^{-1}$ were compared.

Using \eqref{eq:be:admissible}, each bed elevation $b_i$ defines a constraint set $\cK_h^i \subset \cX_h$.  A time-dependent, Stokes-dynamics simulation of $T=200$ years was run for each combination of set $\cK_h^i$ and value $a_j$.  One thousand sample pairs were taken from all runs with a given $\cK_h^i$, at random times in $[0,T]$.  The different $a_j$ values caused advance or retreat of the ice margins; for $a_1<0$ the glacier disappears during the simulation.  Thus the samples included very-different surface elevation pairs.  The initial surface for these runs is a Halfar profile\footnote{For the $b_1=a_2=0$ case the exact time-dependent surface elevation solution \emph{under shallow dynamics} would be known.  Though Stokes dynamics was used in the numerical experiment, the computed surfaces agree closely with this shallow exact solution, as expected; compare \cite{LofgrenAhlkronaHelanow2022}.} \cite{Halfar1981} (Figure \ref{fig:cases}, dotted).  Each time-step was semi-implicit, approximating VI \eqref{eq:be:vi} using $\bu|_{s^{n-1}}\cdot\bn_{s^n}$, with the free-surface stabilization technique from \cite{LofgrenAhlkronaHelanow2022} to lengthen time steps.  (This scheme is only conditionally-stable.)  The VI problems at each step were solved using a reduced-space Newton method with line search \cite{BensonMunson2006}.  The Stokes problems \eqref{eq:glenstokes:weak}, with $\mu_0=10^{-19}\, \text{s}^{-2}$ in \eqref{eq:glen}, were solved using vertically-extruded quadrilateral meshes (Figure \ref{fig:fe:domain}) and stable $Q_2\times Q_1$ Taylor-Hood elements \cite{Elmanetal2014}.

\begin{figure}[ht]
\centering
\includegraphics[width=0.40\textwidth]{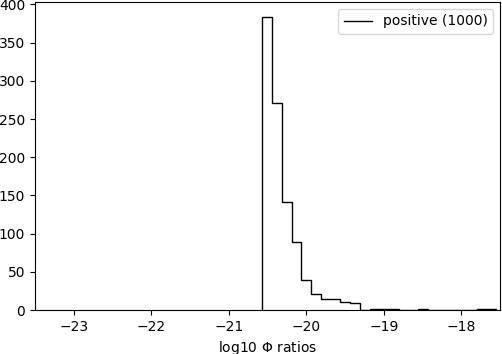} \, \includegraphics[width=0.40\textwidth]{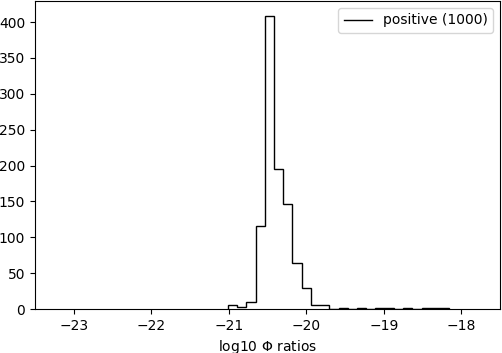}
\caption{Histograms of ratios $\rho^{0.01}(\sigma,s)$ from 1000 sample pairs, for flat-bed (left) and rough-bed (right) cases.}
\label{fig:regratios}
\end{figure}

The essential result is shown in Figure \ref{fig:regratios}, with histograms of sample ratios $\rho^{0.01}(\sigma,s)$, on a logarithmic axis, from the highest spatial resolution ($\Delta x=500$ m and 40 elements in each column).  All ratios are positive, and a coercivity constant $\alpha=10^{-22}$ is credible.  Note that dimensional values of this approximate magnitude are expected because $\|s\|_{\cX}\sim 10^{5}$ is typical in $\cX = W^{1,4}(\Omega)$, using the norm scaling \eqref{eq:norm:Omega}, while residuals in the numerator of \eqref{eq:Phiratio} have typical magnitudes between $10^{-5}$ and $10^0$.  Regarding Conjecture \ref{conj:lipschitz}, Lipschitz continuity for the surface velocity trace, the ratios $\big\|\bu|_\sigma - \bu|_s\big\|_{L^{4/3}}/\|\sigma-s\|_{\cX}$ never exceeded $10^{-6}$.

\begin{figure}[ht]
\centering
\includegraphics[width=0.40\textwidth]{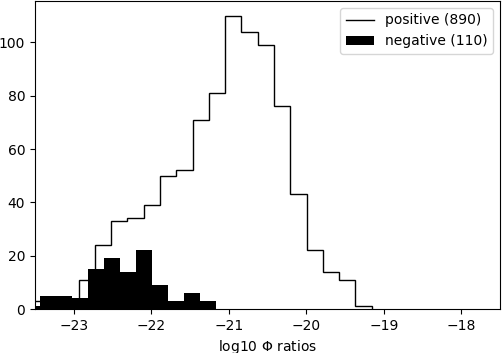} \, \includegraphics[width=0.40\textwidth]{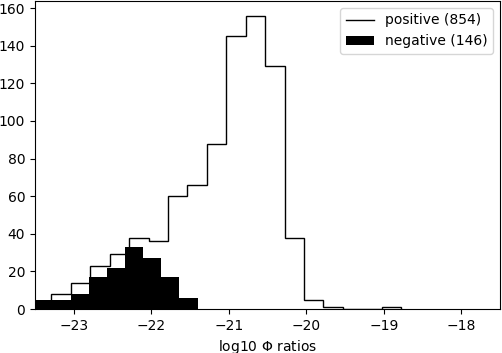}
\caption{Unregularized ratios $|\rho^0(\sigma,s)|$ for the same sample pairs as in Figure \ref{fig:regratios}.}
\label{fig:noregratios}
\end{figure}

In summary we have substantial, but necessarily incomplete, evidence for coercivity Conjecture \ref{conj:regcoercive}.  An underlying reason is that \emph{the regularization does not need to do much work}.  Figure \ref{fig:noregratios} shows unregularized ratios from the same pairs.  About $10\%$ of ratios are negative, which violates coercivity, but their absolute values are an order of magnitude smaller than the positive ratios.  Glacier evolution under Stokes dynamics seemingly ``wants'' to be coercive.

\begin{figure}[ht]
\mbox{\includegraphics[width=0.31\textwidth]{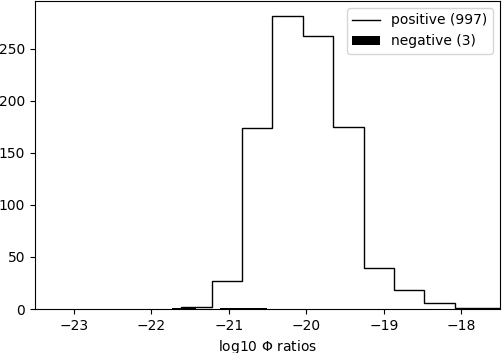} \, \includegraphics[width=0.31\textwidth]{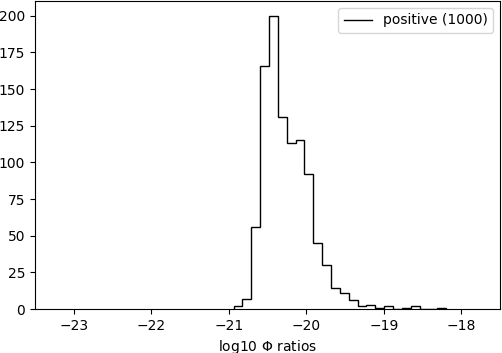} \, \includegraphics[width=0.31\textwidth]{figs/brough500mREG.png}}
\caption{Resolution effects: Ratios $|\rho^{0.01}(\sigma,s)|$ for horizontal mesh resolutions $\Delta x=2000$ m (left), $1000$ m (middle), and $500$ m (right; same as Figure \ref{fig:regratios} right).}
\label{fig:resolutionratios}
\end{figure}

Furthermore the continuum problem appears to be more coercive than any discretization.  Figure \ref{fig:resolutionratios} shows the results from running the same rough-bed experiment as in Figure \ref{fig:regratios}, but at lower resolutions; the coarsest discretization produced a few negative ratios.  Higher marginal resolution, in particular, improves numerical coercivity.

\end{document}